\documentclass[11pt]{article}
\usepackage{amsmath}
\usepackage{amscd}
\usepackage{amsfonts}
\usepackage{amssymb}
\usepackage{amsthm}

\newcommand{\rar}{\rightarrow}
\newcommand{\calg}{\mathcal}

\newcommand{\strc}{\mathcal{O}_{X}}
\newcommand{\strcc}[1]{\mathcal{O}_{#1}}

\newcommand{\enn}{\text{End}}

\newcommand{\dif}{{\mathcal {D}}\text{iff}}
\newcommand{\cale}{\mathcal{E}}
\newcommand{\difdt}{{{\mathcal {D}}\text{iff}}^{\bullet}}
\newcommand{\prdf}{\mathcal{D}}
\newcommand{\diff}{\text{Diff}}
\newcommand{\diffdt}{{\text{Diff}}^{\bullet}}
\newcommand{\hh}[2]{\text{H}^{#1}(X,#2)}
\newcommand{\hoch}[1]{\text{HH}_{#1}}
\newcommand{\choch}[1]{\widehat{\text{HH}}_{#1}}
\newcommand{\ieee}{{\calg I}_{\cale}}

\newcommand{\hcc}[2]{\text{C}^{#1}(#2)}

\newcommand{\cyc}[1]{\text{HC}_{#1}}
\newcommand{\ccyc}[1]{\widetilde{\text{HC}}_{#1}}
\newcommand{\ccccyc}[1]{\widehat{\text{HC}}_{#1}}
\newcommand{\compl}{\mathbb C}
\newcommand{\dol}[2]{{\text{K}^{\bullet}}_{{#1 \text{ } #2}}}
\newcommand{\hhhh}[2]{\text{H}^{#1}(#2)}
\newcommand{\holl}[1]{\text{H}^{#1}}
\newcommand{\ainf}{\text{A}_{\infty}}
\newcommand{\htens}{\widehat{\otimes}}
\newcommand{\id}{\text{id}}
\newcommand{\hcoh}[2]{{\mathbb H}^{#1}(X, #2)}
\newcommand{\fls}{I_{\text{FLS}}}

\newtheorem{thm}{Theorem}
\newtheorem{prop}{Proposition}
\newtheorem{cor}{Corollary}
\newtheorem{lem}{Lemma}

\title{ Some notes on the Feigin-Losev-Shoikhet integral conjecture}
\author{Ajay C. Ramadoss}

\begin{document}

\maketitle

\begin{abstract}
Given a holomorphic vector bundle $\cale$ on a smooth connected
compact complex manifold $X$, [FLS] use a notion of completed
Hochschild homology $\choch{}{}$  of $\dif(\cale)$ such that
$\choch{0}{(\dif(\cale))}$ is isomorphic to \\ $\hh{2n}{\compl}$. On
the other hand, they construct a linear functional on
$\choch{0}{(\dif(\cale))}$. This therefore gives rise to a linear
functional $I_{\cale}$ on $\hh{2n}{\compl}$. They show that this
functional is $\int_X$ if $\cale$ has non zero Euler characteristic.
They conjecture that this functional is $\int_X$ for all $\cale$. \\

These notes prove that $I_{\cale} = I_{\calg F}$ for any pair $(\cale, \calg F)$ of holomorphic vector bundles on $X$.
In particular, if $X$ has one vector bundle with nonzero Euler characteristic, $I_{\cale} = \int_X$ for every vector bundle $\cale$ on $X$.
\\

[FLS] also use a notion of completed cyclic homology $\ccccyc{}{}$
of $\dif(\cale)$ such that $\ccccyc{-i}{(\dif(\cale))} \simeq
\hh{2n-i}{\compl} \oplus \hh{2n-i+2}{\compl} \oplus .....$ The
construction yielding $I_{\cale}$ generalizes to yield linear
functionals on \\ $\ccccyc{-2i}(\dif(\cale))$ for each $i \geq 0$.
The linear functional thus obtained on $\ccccyc{-2i}(\dif(\cale))$
thus yields a linear functional $I_{\cale, 2i,2k}$ on
$\hh{2n-2k}{\compl}$ for $0 \leq k \leq i$. [FLS] conjecture that
$I_{\cale,2,0} = \int_X$
and make a further conjecture about $I_{\cale,2,2}$. \\

These notes prove that $I_{\cale,2i,0} = I_{\cale}$ for all $i \geq
0$. In particular, if $X$ has at least one vector bundle with non
zero Euler characteristic, then $I_{\cale,2i,0} = \int_X$. We also
prove that $I_{\cale,2i,2k} = 0$ for $k
> 0$. The latter is stronger than what [FLS] expect when $i=k=1$.  \\



\end{abstract}

\section{Introduction}

 Let $X$ be a smooth connected compact complex manifold such that $\text{dim}_{\mathbb C} X = n$. Let ${\calg E}$ be a holomorphic vector bundle on $X$. Let $\dif(\calg E)$ be the
 sheaf of {\it holomorphic} differential operators on $\calg {E}$. Let $\diff(\cale)$
 be the algebra of global sections of $\dif(\cale)$. One is interested in
understanding the Hochschild homology of $\diff(\cale)$. This is
however, not an object that one can easily analyze. \\

Let $\text{Dolb}^{\bullet}(X,\strc)$ denote the Dolbeaux resolution
of the sheaf $\strc$ of holmorphic functions on $X$. Let
$\dif^{\bullet}(\cale) = \text{Dolb}^{\bullet}(X,\strc)
\otimes_{\strc} \dif(\cale)$. Following [FLS] we replace
$\diff(\cale)$ by
 $$ \diffdt(\cale) := \Gamma(X,\dif^{\bullet}(\cale) ) \text{ . } $$

Let $\dol{\cale}{ }$ denote the complex
$\Gamma_{\text{C}^{\infty}}(X,\text{Dolb}^{\bullet}(X,\strc)
\otimes_{\strc} \cale)$ of smooth global sections of
$\text{Dolb}^{\bullet}(X,\strc) \otimes_{\strc} \cale$ . This is the
Dolbeaux complex of $\cale$. If $\calg D$ is a holomorphic
differential operator on $\cale$, then $\calg D$ gives rise to a
$\compl$-endomorphism of the complex $\dol{\cale}{ }$. It follows
that $\calg D$ induces an endomorphism ${\calg D}_*: \hh{i}{\cale}
\rar \hh{i}{\cale}$ for each $i$. The {\it supertrace} of $\calg D$
is given by
$$\text{str}({\calg D}) := \sum_i {(-1)}^{i}
 Tr_{\hh{i}{\cale}} {\calg D}_*  \text{ . }$$

 Note that $\diff(\cale)$ is a subalgebra of $\diffdt(\cale)$. It follows that there is a natural
 map from the $i$th Hochschild homology of $\diff(\cale)$ to that of
$\diffdt(\cale)$ for all $i$. Also recall that if $\calg A$ is any
associative $\compl$-algebra concentrated in degree $0$ then the
$0$-th Hochschild homology $\hoch{0}(\calg A)$ of $\calg A$ is the
Abelianization of $\calg A$.
 It follows that any trace on $\calg A$ is a linear functional on
$\hoch{0}(\calg A)$ and a Hochschild $0$-cocycle of $A$ .Note that
the supertrace described in the previous paragraph
is a Hochschild $0$-cocycle of $\diff(\cale)$.\\

 In [FLS], a Hochschild $0$-cocycle $\text{tr}$ on $\diffdt(\cale)$ is constructed algebraically.
 This Hochschild cocycle is constructed in such a way that the following diagram
 commutes.  $$\begin{CD}
 \hoch{0}(\diff(\cale)) @>>> \hoch{0}(\diffdt(\cale)) \\
 @V{\text{str}}VV   @VV{\text{tr}}V \\
 \compl  @>>{\text{id}}> \compl \\
 \end{CD} $$
The top row of the above diagram is the natural map from
$\hoch{0}(\diff(\cale))$ to $\hoch{0}(\diffdt(\cale))$ arising out
of the natural inclusion of $\diff(\cale)$ as a subalgebra of
$\diffdt(\cale)$. \\

The theory that is more related to the cohomology of $X$ is however,
that of the {\it completed Hochschild homology} of the sheaf
$\dif(\cale)$ of holomorphic differential operators on $\cale$ .
This is defined in this paper in Section 3. Let
$\choch{i}(\dif(\cale))$ denote the $i$ -th completed Hochschild
homology of $\dif(\cale)$. By a lemma (Lemma 3) that generalizes a
Theorem of Brylinski [Bryl],
$$\choch{-i}(\dif(\cale)) \simeq \hh{2n-i}{\compl} \text{ . }$$ Call this
isomorphism $\beta_{\cale}$.This is the key property of
$\choch{*}(\dif(\cale))$ that relates it to the cohomology of $X$
with coefficients in $\compl$. [FLS] also use a notion of completed
Hochschild homology that has the same relation with
$\hh{*}{\compl}$. However, as hinted in [FLS] itself, understanding
the isomorphism $\beta_{\cale}$ is easier with the definition of
completed Hochschild homology used
here. \\

[FLS] shows that the Hochschild cocycle $\text{tr}$ of
$\diffdt(\cale)$ "extends" to yield a linear functional
$\hat{\text{tr}}$ on $\choch{0}(\dif(\cale))$. In this paper, we are
forced to redo this very step in Section 3 since we use a definition
of $\choch{*}(\dif(\cale))$ that is {\it a priori } different from
that of [FLS]. Our definition of $\choch{*}(\dif(\cale))$ implies
that there is a natural map from $\hoch{0}(\diff(\cale))$ to
$\choch{0}(\dif(\cale))$ (see Section 3 for more details). The
outcome of our "extending" $\text{tr}$ to a linear functional
$\hat{\text{tr}}$ on
$\choch{0}(\dif(\cale))$ is the following commutative diagram.\\

$$\begin{CD}
\hoch{0}(\diff(\cale)) @>>> \choch{0}(\dif(\cale))\\
@VV{\text{str}}V    @V{\hat{\text{tr}}}VV \\
\compl @>>{\text{id}}> \compl \\
\end{CD}$$

The linear functional we construct on $\choch{0}(\dif(\cale))$ still
uses the basic construction of [FLS] and , as will be clear later in
this paper, yields an algebraic construction of $\int_X$ when $X$
has
at least one vector bundle with non zero Euler characteristic. \\

Denote the composite map $\hat{\text{tr}} \circ \beta_{\cale}^{-1}:
\hh{2n}{\compl} \rar \compl$ by $I_{\cale}$. $I_{\cale}$ is a linear
functional on $\hh{2n}{\compl}$. \\

 The following theorem was proven in [FLS]. Their proof goes through
 in our setup as well.\\

 \begin{thm} [FLS] If the Euler characteristic $\chi_{\cale}$ of $\cale$ is nonzero, then $I_{\cale} =
 \int_X$.
 \end{thm}

Further, in [FLS], it was conjectured that $I_{\cale} = \int_X$ for
all vector bundles $\cale$ , even those of zero Euler
characteristic. In these notes, we prove the following theorem.

\begin{thm} If $\cale$ and $\mathcal G$ are two vector bundles on $X$, then $I_{\cale} = I_{\mathcal
G}$.
\end{thm}

The following corollary is immediate from this and Theorem 1.

\begin{cor} If $X$ has a vector bundle $\mathcal G$ whose Euler characteristic $\chi_{\mathcal G}$ is nonzero, then for any vector bundle
$\cale$ on $X$, $I_{\cale} = \int_X$.
\end{cor}

In particular, if $X$ is a smooth separated scheme over $\compl$,
then the skyscraper sheaf supported at a point $p$ of $X$ has
nonzero Euler characteristic. On the other hand, it is a coherent
$\strc$-module. It therefore, has a finite resolution by
(holomorphic) vector bundles on $X$. It follows that at least one
vector bundle in such a resolution has non-zero Euler
characteristic. Thus, if  $\cale$ is a vector bundle on a smooth
compact complex manifold $X$ that is a smooth separated scheme over
$\compl$, then $I_{\cale}= \int_X$.\\

Note that if $\calg D$ is a global holomorphic differential operator
on a vector bundle $\cale$ on $X$, then $\calg D$ is a $0$-cycle of
the Hochschild chain complex of $\diff(\cale)$. It follows that
$\calg D$ yields an element in $\choch{0}(\dif(\cale))$, and hence,
an element in $\hh{2n}{\compl}$ which we will denote by $[\calg D]$.
It follows from Corollary 1 that if $X$ admits at least one vector
bundle with non-zero Euler characteristic, then
$$\text{str}(\calg D) = \int_X [\calg D] \text{ . } $$

The above statement is somewhat similar to a result (Corollary 5.6)
in a paper by P.Schapira and J-P.Schneiders [S-S]. More recently, an
entirely different approach has yielded a proof of a stronger
version of this statement (the supertrace theorem)without the
condition that $X$ admit a vector bundle with non-zero Euler
characteristic [EnFe]. In a sequel to this paper, we shall develop
our approach further to show that $I_{\cale} = \int_X$ for any
vector bundle $\cale$ on any compact complex manifold $X$. This
shall yield another proof of the supertrace theorem different from
the one in [EnFe]. \\

One also defines the completed cyclic homology
$\ccccyc{*}(\dif(\cale))$ of $\dif(\cale)$. Further, as in [FLS],
the construction of the linear functional on the $0$-th completed
Hochschild homology of $\dif(\cale)$ is extended to the construction
of
multiple linear functionals $\hat{\text{tr}}_{2i} : \ccccyc{-2i}(\dif(\cale)) \rar \compl$. On the other hand, it is shown that one has an isomorphism \\
$$ \ccccyc{-2i}(\dif(\cale)) \simeq \hh{2n}{\compl} \oplus \hh{2n-2}{\compl} \oplus ..... \oplus \hh{2n-2i}{\compl} \text{ . }$$

Denote the inverse of this isomorphism by $\ieee^i$. Let
$I_{\cale,2i} = \hat{\text{tr}_{2i}} \circ \ieee^i$. Then,
$$ I_{\cale,2i} :  \hh{2n}{\compl} \oplus \hh{2n-2}{\compl} \oplus .... \oplus \hh{2n-2i}{\compl} \rar \compl \text{ . } $$

For $0 \leq k \leq i$ , let $I_{\cale, 2i, 2k}: \hh{2n-2k}{\compl}
\rar \compl $ be the composition of $I_{\cale, 2i}$ with the
inclusion of $\hh{2n-2k}{\compl}$ in $\hh{2n}{\compl} \oplus
\hh{2n-2}{\compl} \oplus ..... \oplus \hh{2n-2i}{\compl}$ as a
direct summand . In [FLS] it was conjectured that $I_{\cale, 2, 0} =
C. \int_X$ where $C$ is a constant independent of $\cale$. Further,
there was a prediction in [FLS] that $I_{\cale,2,2}$ is the integral
over a cycle related to the Chern classes of $\cale$ and those of
vector bundles intrinsic to $X$. In this note,
we prove the following theorem. \\

\begin{thm} $I_{\cale, 2i, 0} =  I_{\cale}$ . Further, $I_{\cale, 2i, 2k} = 0$ if
$k>0$.
\end{thm}

The latter part is more than what [FLS] expected in their work,
though a piece of circumstantial evidence for $I_{\cale,2,2}$ to be
$0$ is provided there. By this theorem and Corollary 1, $I_{\cale,
2i, 0} = \int_X$ whenever $X$ has a vector bundle with nonzero Euler
characteristic. In particular,
this happens when $X$ is a smooth separated scheme over $\compl$. \\

\subsection{Outline of this note} The proofs of theorems 2 and 3 are simple and involve two easy "bookkeeping" lemmas (Lemma 1  and Lemma 4 )
. These are done in Sections 3 and 4 respectively.  \\

Section 2 recalls basic notions about Hochschild homology that we
require. In addition, Section 2 recalls the definition of an
$\ainf$-morphism between two differential graded algebras. It shows
that an $\ainf$-morphism $\calg F$ between two dg-algebras $\calg A$
and $\calg B$ induces a map ${\calg F}_{\text{hoch}}$ of complexes
from the complex of Hochschild chains of $\calg A$ to that of $\calg
B$ (Proposition 2 ). No material in this section is new. Proposition
2 is from [FLS] .\\

The first part of Section 3 recalls the basic set up for this note
from [FLS]. In particular, we recall the construction of the
Hochschild $0$-cocycle $\text{tr}$ of $\diffdt(\cale)$ from [FLS].
The completed Hochschild homology $\choch{*}(\dif(\cale))$ is
defined in Section 3. We show that the formula for the Hochschild
cocycle $\text{tr}$ "extends" to a formula for a linear functional
$\hat{\text{tr}}$ on $\choch{0}(\dif(\cale))$. [FLS] in [1] do
something very similar. The definition of $\choch{*}(\dif(\cale))$
that is used in this paper is however, a priori somewhat different
from that used in [FLS]. This forces us to redo the parts where they
extend the Hochschild cocycle $\text{tr}$ of $\diffdt(\cale)$ to a
linear functional
on $\choch{0}(\dif(\cale))$.  \\

The last part of Section 3 is devoted to the "bookkeeping Lemma
for $\hat{\text{tr}}$" (Lemma 1). \\

Section 4 is devoted to understanding the relation between
$\choch{*}(\dif(\cale))$ and $\hh{*}{\compl}$. In particular, we
prove that $$\choch{-i}(\dif(\cale)) \simeq \hh{2n-i}{\compl}$$
(Corollary 6 ). For $\cale = \strc$, this was proven by Brylinski
[Bryl]. Let $\beta_{\cale}$ denote the isomorphism between
$\choch{-i}(\dif(\cale))$ and $\hh{2n-i}{\compl}$. Section 4 also
proves the "bookkeeping Lemma for $\beta$" (Lemma 4 ). Theorem 2  is
an
immediate consequence of Lemma 1 and Lemma 4 . \\

 Section 5 has four subsections. The first subsection recalls the basic notions of cyclic
 homology that we require. The second one carefully examines the linear functionals on the completed
 cyclic homology of $\dif(\cale)$. The third subsection examines the relation between the completed
 cyclic homology of $\dif(\cale)$ and $\hh{*}{\compl}$. The final subsection writes down the final steps that prove Theorem 3. \\

 \textbf{Note:} All complexes that appear in this paper are cochain
 complexes by convention. The term {\it Hochschild cocycle } of a dg-algebra shall
  only refer to cocycles of the complex of Hochschild cochains of
 that algebra as defined in Definition 2, Section 2.1.\\

 \textbf{Acknowledgements:} This paper and its sequel would not have been what they are but for the help I received from many quarters.
 I am very grateful to Prof. Boris
Tsygan for going through the paper carefully, introducing me to a
paper of N. Teleman [Tel] and for some very useful comments and
suggestions. I am also very grateful to Prof. Madhav Nori for some
very useful discussions, comments and suggestions. Heartfelt thanks
are also due to Prof. Ryszard Nest and Prof. Alexander Gorokhovsky
for helping me understand the correct completed tensor product to be
used.
 I am also grateful to Prof. Shrawan Kumar and Dr. Victor Protsak for useful discussions and to the referee for going through this paper carefully
 and helping me streamline its presentation. \\

\section{Basic results about Hochschild homology}

{\it Conventions used in this section}\\

 1. The term "dg-algebra" in
this section refers to a differential
graded $\compl$-algebra with unit. \\

2. All complexes of $\compl$-vector space are cochain complexes,
i.e, the differential of any complex of $\compl$-vector spaces has
degree $+1$. \\

\subsection{The Hochschild chain complex and Bar complex of a
dg-algebra}

{\it Definition 1: } If $\calg A$ is a differential graded $\compl$-algebra with differential $\delta$, the complex of Hochschild chains
$\hcc{\bullet}{\calg A}$ is the cochain complex obtained by equipping the graded $\compl$-vector space $\oplus_{i \geq 1} {\calg A [1]}^{\otimes
i} [-1]$ with the {\it Hochschild differential}. The Hochschild differential $d$ is given by the formula $$ d(a_0 \otimes .... \otimes a_n) =
\sum_{i=0}^{i=n-1} {(-1)}^{(d_0 + ... +d_i + i+1)} a_0 \otimes ... \otimes a_ia_{i+1} \otimes .... \otimes a_n $$ $$ +
{(-1)}^{(d_n+1)(d_0+...+d_{n-1}+n-1)} a_na_0 \otimes a_1 \otimes ... \otimes a_{n-1} $$ $$+ \sum_{j=0}^{j=n} {(-1)}^{(d_0 + ... +d_{j-1}+j)} a_0
\otimes ... \otimes \delta(a_j) \otimes ... \otimes a_n $$
  for homogenous elements $a_0,...,a_n$ of $\calg A$ of degrees $d_0,..,d_n$ respectively. \\

  Let $\tau_n:{\calg A}^{\otimes n} \rar {\calg A}^{\otimes n}$
  denote the endomorphism given by the formula $$a_1 \otimes .... \otimes a_n \leadsto
  {(-1)}^{(d_n+1)(d_1+...+d_{n-1}+n-1)}a_n \otimes a_1 \otimes .... \otimes
  a _{n-1} $$ for homogenous elements $a_1,...,a_n$ of $\calg A$ of degrees $d_1,..,d_n$ respectively.
   Let $\partial_1:{\calg A}^{\otimes n} \rar {\calg A}^{\otimes n-1}$ denote the
   morphism given by the formula
  $$a_1 \otimes .... \otimes a_n \leadsto {(-1)}^{d_1+1}a_1a_2 \otimes ... \otimes a_n $$ for
   homogenous elements $a_1,...,a_n$ of $\calg A$ of degrees $d_1,..,d_n$ respectively.Then,
   the Hochschild differential on ${\calg A}^{\otimes k}$ is also given by the
  formula $$\sum_{i=1}^{i=k}
\tau_{k-1}^{i-1} \circ \partial_1 \circ \tau_k^{k-i+1} +
\sum_{i=1}^{i=k} \tau_{k}^{i-1} \circ (\delta \otimes id \otimes..
\otimes id) \circ \tau_k^{k-i+1} \text{ . }$$

Contrary to the standard practice, we will refer to a degree $-n$
cocycle of $\hcc{\bullet}{\calg A}$ as a {\it Hochschild $n$- cycle
of $\calg
A$} .\\

 The Hochschild homology $\hoch{-n}(\calg A)$ of $\calg A$ is the $-n$th cohomology of the complex $\hcc{
\bullet}{\calg A}$.  By the definition of the Hochschild complex, it
is clear that if $\calg A$ is a $\compl$-algebra (i.e, $\calg A$ is
concentrated in degree $0$), then $\hoch{0}(\calg A) = \frac{\calg
A}{[\calg A, \calg A]}$ where $[\calg A, \calg A]$ is the commutator
of $\calg A$.Also, if $\calg A$ is a $\compl$-algebra,
$\hoch{-n}(\calg A) = 0$ if $n<0$. Also, if $\calg A$ is dg-algebra,
then elements of $\text{Ker}(\delta) \subset {\calg A}^0$ are
$0$-Hochschild cycles of $\hcc{\bullet}{\calg A}$ where $\delta$
denotes the internal differential on $\calg A$. It follows that there is a canonical map $\text{Ker}(\delta) \rar \hoch{0}(\calg A)$. \\

 {\it Definition 2: } The {\it complex of Hochschild cochains } of $\calg A$ is the complex whose $-i$-th term is
 $\text{Hom}_{\compl}(\hcc{i}{\calg A},\compl)$ with differential induced by that on $\hcc{\bullet}{\calg A}$. A {\it Hochschild cocycle}
 of ${\calg A}$ is a cocycle of the complex of Hochschild cochains of $\calg A$. Note that a Hochschild $0$-cocycle of $\calg A$ induces a
 $\compl$-linear functional on $\hoch{0}(\calg A)$.   \\

{\it Definition 3:} The {\it Bar complex}
$\text{bar}^{\bullet}(\calg A)$ of $\calg A$ is the cochain complex
obtained by equipping the graded $\compl$-vector space $\oplus_{i
\geq 1} {\calg A [1]}^{\otimes i}[-1]$ with the {\it Bar
differential}. The Bar differential is given by the formula
$$ d(a_0 \otimes .... \otimes a_n) =  \sum_{i=0}^{i=n-1}
{(-1)}^{(d_0 + ... +d_i + i+1)} a_0 \otimes ... \otimes a_ia_{i+1} \otimes .... \otimes a_n $$ $$+ \sum_{j=0}^{j=n} {(-1)}^{(d_0 + ...
+d_{j-1}+j)} a_0 \otimes ... \otimes \delta(a_j) \otimes ... \otimes a_n $$  for homogenous elements
$a_0,...,a_n$ of $\calg A$ of degrees $d_0,..,d_n$ respectively. \\

 The Bar differential on ${\calg A}^{\otimes k}$ is also given by the
  formula $$\sum_{i=1}^{i=k-1}
\tau_{k-1}^{i-1} \circ \partial_1 \circ \tau_k^{k-i+1}  +
\sum_{i=1}^{i=k} \tau_{k}^{i-1} \circ (\delta \otimes id \otimes..
\otimes id) \circ \tau_k^{k-i+1} \text{ . }$$

 {\it Recollection 1 :} If $s_n: \text{bar}^{\bullet}(\calg A) \rar
\text{bar}^{\bullet-1}(\calg A)$ is the map $a_0 \otimes ... \otimes
a_{n} \leadsto -1 \otimes a_0 \otimes ... \otimes a_{n}$ then the
maps $\{s_n\} \text{ , } n \geq -1$ give a homotopy between $id:
\text{bar}^{\bullet}(\calg A) \rar \text{bar}^{\bullet}(\calg A)$
and $0: \text{bar}^{\bullet}(\calg A) \rar
\text{bar}^{\bullet}(\calg A)$. It follows that
$\text{bar}^{\bullet}(\calg A)$ is acyclic. \\

If $V^{\bullet}$ is a finite dimensional graded $\compl$-vector
space then \\ $\enn(V^{\bullet}) = \oplus_{i,j}
\text{Hom}(V^i,V^j)$. Let $\pi_{i,j}$ denote the projection from
$\enn(V^{\bullet})$ to $\text{Hom}(V^i,V^j)$. If $M \in
\enn(V^{\bullet})$ , the supertrace of $M$ is
the alternating sum $\sum_i {(-1)}^i tr(\pi_{i,i}(M))$. We recall the following proposition from [FLS]. \\

\begin{prop} Let $V^{\bullet}$ be a finite dimensional graded $\compl$-vector space with zero differential. \\

1. $\hoch{i}(\enn(V^{\bullet})) = 0$ for $i \neq 0$ \\
2. $\hoch{0}(\enn(V^{\bullet})) \simeq \compl$ \\
3. The composite $\enn(V^{\bullet})^0 \rar \hoch{0}(\enn(V^{\bullet})) \rar \compl$ takes
an element of $\enn(V^{\bullet})^0$ to its supertrace. \\
4. The isomorphism $\hoch{0}(\enn(V^{\bullet})) \simeq \compl$ takes
the class in $\hoch{0}(\enn(V^{\bullet}))$ of a Hochschild $0$-cycle
in $\enn(V^{\bullet})^{\otimes k}$ to $0$ for all $k \geq 2$.

\end{prop}

\begin{proof} The proof of this proposition is a trivial modification of the proof of the Morita invariance of
Hochschild homology for matrices in Loday [2](see Theorem 1.2.4 of [2]).
\end{proof}

We may therefore, denote the isomorphism $\hoch{0}(\enn(V^{\bullet})) \simeq \compl$ by $\text{str}$. \\

\subsection{$\ainf$-morphisms between dg-algebras}

Let $\calg A$ and $\calg B$ be two dg-algebras. Let ${\calg
F}_i:{\calg A}^{\otimes i} \rar {\calg B}[n_i]$ be $\compl$-linear
maps with $n_i \in \mathbb Z$ for each $i$. If $k_1,...,k_l$ are
positive integers such that $\sum_j k_j = k$ and if $\sum_j n_j =n$,
then ${\calg F}_{k_1} \boxtimes ... \boxtimes {\calg F}_{k_l}$ will
denote the $\compl$-linear map from ${\calg A}^{\otimes k}$ to
${\calg B}^{\otimes l}[n]$ such that $$ {\calg F}_{k_1} \boxtimes
... \boxtimes {\calg F}_{k_l}(a_1 \otimes .. \otimes a_k) = {\calg
F}_{k_1}(a_1 \otimes ... \otimes a_{k_1}) \otimes .... \otimes
{\calg F}_{k_l}(a_{k_1+.....k_{l-1}+1} \otimes ... \otimes a_{k})
\text{ . }$$

{\it Definition 2 :} An $\ainf$-morphism $\calg F$ from a dg-algebra
$\calg A$ to a dg-algebra $\calg B$ is a collection of maps $${\calg
F}_k : {\calg A}^{\otimes k} \rar {\calg B}[1-k] $$ for all $k \geq
1$ such that the map ${\calg F}_{\text{bar}}:
\text{bar}^{\bullet}(\calg A) \rar \text{bar}^{\bullet}(\calg B) $
defined by $${\calg F}_{\text{bar}}(a_1 \otimes ... \otimes a_k) =
\sum_{ \{(k_1,...,k_l) \text{ } | l >0 \text{ and } \sum_j k_j =k \}
} {\calg F}_{k_1} \boxtimes ... \boxtimes {\calg F}_{k_l}(a_1
\otimes ... \otimes a_k) $$ is a morphism of complexes from
$\text{bar}^{\bullet}(\calg A)$ to $\text{bar}^{\bullet}(\calg B)$. The maps ${\calg F}_k$ are called the {\it Taylor components} of ${\calg F}$. \\

The condition that ${\calg F}_{\text{bar}}$ commutes with the
differentials on the bar complexes of ${\calg A}$ and ${\calg B}$
respectively is equivalent to the condition that the maps ${\calg
F}_k$ satisfy the following relations.

$$ [\pm {\calg F}_{k-1}(a_1.a_2 \otimes a_3 \otimes .... \otimes a_k) \mp {\calg F}_{k-1}(a_1 \otimes a_2. a_3 \otimes .... \otimes a_k) \pm $$
$$ \pm {\calg F}_{k-1}(a_1 \otimes a_3 \otimes .... \otimes a_{k-1}.a_k) \pm {\calg F}_{k}(\delta(a_1 \otimes a_2 \otimes .... \otimes a_k)) \pm
$$
$$\pm \delta({\calg F}_k(a_1 \otimes ... \otimes a_k)) \pm \sum_{l=1}^{l=k-1} \pm {\calg F}_{l}(a_1 \otimes .... \otimes a_l) \circ {\calg F}_{k-l}(a_{l+1} \otimes .... \otimes a_k) = 0$$
$$k \geq 1 $$

We recall the following proposition from [FLS]. \\

\begin{prop} An $\ainf$ morphism $\calg F$ from an associative dg-algebra $\calg A$ to an associative dg-algebra $\calg B$ induces a map ${\calg
F}_{\text{Hoch}}$ of complexes from the $\hcc{\bullet}{\calg A}$ to $\hcc{\bullet}{\calg B}$.
\end{prop}

We recall the proof of this proposition from [FLS].\\

\begin{proof}

Let $\tau:{\calg A}^{\otimes k} \rar {\calg A}^{\otimes k}$ be the map which takes $a_1 \otimes ... \otimes a_k$ to
 \\ ${(-1)}^{(d_k+1)(d_1+...+d_{k-1}+k-1)}a_k \otimes a_1 \otimes ... \otimes a_{k-1}$ for homogenous elements $a_1,...,a_k$ of $\calg A$ of degrees
$d_1,...,d_k$ respectively. Consider the map ${\calg F}_{\text{Hoch}}:\hcc{\bullet}{\calg A} \rar \hcc{\bullet}{\calg B}$ defined by the formula
$$ {\calg F}_{\text{Hoch}}(a_1 \otimes .... \otimes a_k) = \sum_{
\{(k_1,...,k_l) \text{ } | l >0 \text{ and } \sum_j k_j =k \} }
[{\calg F}_{k_1} \boxtimes ... \boxtimes {\calg F}_{k_l}(a_1 \otimes
... \otimes a_k) $$ $$ + \sum_{j=1}^{j=k_l-1} {\calg F}_{k_l}
\boxtimes {\calg F}_{k_1} \boxtimes ... \boxtimes {\calg
F}_{k_{l-1}}(\tau^j(a_1 \otimes ... \otimes a_k)) ] \text{ . }$$

We leave the verification that ${\calg F}_{\text{Hoch}}$ respects the Hochschild differential to the reader. \\

\end{proof}

Let $\calg A$ be a dg-algebra. Let $V^{\bullet}$ be a finite dimensional graded $\compl$-vector space with $0$ differential. Suppose that $\calg
F$ is an $\ainf$-morphism from $\calg A$ to $\enn(V^{\bullet})$ with Taylor components ${\calg F}_k$. Recall that for any $k >0$, we have a map
$\tau:{\calg A}^{\otimes k} \rar {\calg A}^{\otimes k}$ such that $\tau(a_1 \otimes ... \otimes a_k) = {(-1)}^{(d_k+1)(d_1+...+d_{k-1}+k-1)}a_k
\otimes a_1 \otimes ... \otimes a_{k-1}$ for homogenous elements $a_1,...,a_k$ of $\calg A$ of degrees $d_1,...,d_k$ respectively. By Proposition 1,
$\text{str}$ is a Hochschild $0$-cocycle of $\enn(V^{\bullet})$. We now have the following Corollary of Proposition 2. \\

\begin{cor} The supertrace on $\enn(V^{\bullet})$ pulls back to a Hochschild $0$-cocycle $\text{tr}$ of $\calg A$. On Hochschild $0$-cycles of $\calg A$ that
arise out of elements on degree $k-1$ in ${\calg A}^{\otimes k}$,
the Hochschild cocycle $\text{tr}$ is given by the map from ${\calg
A}^{\otimes k}$ to $\compl$ given by
 $$a_1 \otimes ... \otimes a_k \leadsto \sum_{j=0}^{j=k-1} \text{str}({\calg F}_k(\tau^j(a_1 \otimes ... \otimes a_k))) \text{ . } $$
  \end{cor}

 \begin{proof} The Hochschild cocycle $\text{tr}$ is given by $\text{tr}(x) = \text{str}({\calg F}_{\text{hoch}}(x))$ for any $x \in \hcc{0}{\calg A}$.
 Note that by Proposition 2, ${\calg F}_{\text{hoch}}(x) \in \hcc{0}{\enn(V^{\bullet})}$. The exact formula for $\text{tr}$ given in this corollary is now immediate from
 the formula for ${\calg F}_{\text{hoch}}$ given in the proof of Proposition 2. \\
 \end{proof}

\section{The completed Hochschild homology of $\dif(\cale)$ and a linear functional on $\choch{0}(\dif(\cale))$.}

\subsection{The basic construction in [FLS].}

Let $\dol{\cale}{ }$ denote the Dolbeaux complex of $\cale$ , as in
the introduction to this paper. Then, $\dol{\cale}{ }$ decomposes as
the direct sum of a complex with zero differential and an acyclic
complex i.e, $\dol{\cale}{} = \dol{0}{\cale} \oplus \dol{1}{\cale}$
where $\hhhh{\bullet}{{\dol{0}{\cale}}} =
\hhhh{\bullet}{\dol{\cale}{ }}$ , $\dol{0}{\cale}$ has $0$
differential and $\dol{1}{\cale}$ is acyclic.
This is a consequence of Hodge theory (for instance,see Theorem 5.24 in [Vois] ). \\

Recall that $\dif^{\bullet}(\cale) = \text{Dolb}^{\bullet}(X,\strc) \otimes_{\strc} \dif(\cale)$ and that
 $$ \diffdt(\cale) := \Gamma(X,\dif^{\bullet}(\cale) ) \text{ . }$$

{\it The main construction in [FLS] is that of an $\ainf$ morphism $\calg F$ from $\diffdt(\cale)$ to $\enn(\dol{0}{\cale})$}. Note that
 $\dol{0}{\cale}$ is a finite dimensional $\compl$-vector space with $0$ differential. As $\compl$-vector spaces, $\text{K}_{0 \cale}^i \simeq
 \text{H}^i(X, \cale)$. We may therefore apply Propositions 1, 2 and Corollary 2 of Section 2 with $\calg A = \diffdt(\cale)$ and $\calg B =
 \enn(\dol{0}{\cale})$. We obtain the following facts immediately. \\

 {\it Fact 1: } By Proposition 1, $\hoch{i}(\enn(\dol{0}{\cale})) = 0 $ for $i \neq 0$ and \\ $\hoch{0}(\enn(\dol{0}{\cale})) \simeq
 \compl$. This isomorphism is induced by the map taking a degree $0$ element of $\enn(\dol{0}{\cale})$ to its supertrace. \\

 {\it Fact 2:} By Proposition 2, the supertrace on $\enn(\dol{0}{\cale})$ pulls back to a Hochschild $0$-cocycle $\text{tr}$ on $\diffdt(\cale)$.
 If $D$ is a Hochschild $0$-cycle
  of $\diffdt(\cale)$ that arises out of a degree $k-1$ element of $\diffdt(\cale)^{\otimes k}$, and if $[D]$ denotes the class of $D$ in
  $\hoch{0}(\diffdt(\cale))$, then $$ \text{tr}([D]) = \sum_{j=0}^{j=k-1} \text{str}({\calg F}_k(\tau^j(D))) \text{ . }$$

  \textbf{Notation:} We shall also denote the map $$ D \mapsto \sum_{j=0}^{j=k-1} \text{str}({\calg F}_k(\tau^j(D)))$$ by $\fls$. \\

\subsubsection{Construction of $\calg F$}

 We now recall the construction of $\calg F$ from [FLS]. Let $C_{k}$ denote the configuration space
$\{ t_1 < ... < t_{k} \text{ } | \text{ } t_i \in {\mathbb R}
\}/G^{(1)}$ where $G^{(1)}$ is the one dimensional group of shifts
$(t_1,...,t_k) \rar (t_1+c,..,t_k+c)$ . This is a smooth $k-1$
dimensional manifold though it is {\it not} compact if $k \geq 2$.
Let $\tau_i = t_{i+1}-t_i$ for $1 \leq i \leq k-1$. The map
$(t_1,...,t_k) \leadsto (\tau_1,...,\tau_{k-1})$ is a diffeomorphism
between $C_k$ and the product
 $\Pi_{i=1}^{i=k-1} \{\tau_i > 0\}$. Let $\overline{\{\tau_i > 0\}}$ denote the  compactification of $\{\tau_i \geq 0\}$
  by a point at infinity. The cube $\overline{C_k} := \Pi_{i=1}^{i=k-1} \overline{\{\tau_i > 0\}} $
  is a compactification of $C_k$. \\

Let $D=D_1 \otimes ... \otimes D_k \in \diffdt(\cale)^{\otimes k}$.
Recall that the $D_i$ yield endomorphisms of $\dol{\cale}{ }$ as
follows: Let $U \subset X$ be an open ball on which $\cale$ is
trivial. Let $z_1,...,z_n$ be local holomorphic coordinates on $U$.
A section of $\dol{\cale}{ } |_U$ is a linear combination of
sections of the form $s \otimes \bar{dz_{i_1}} \wedge ... \wedge
\bar{dz_{i_m}}$.  $D_i(s \otimes \bar{dz_{i_1}} \wedge ... \wedge
\bar{dz_{i_m}}) = D_i(s) \otimes
\bar{dz_{i_1}} \wedge ... \wedge \bar{dz_{i_m}}$. \\

Let ${\bar{\partial_{\cale}}}^*$ be the Hodge adjoint of
$\bar{\partial_{\cale}}$. Let $\Delta_{\cale}$ denote the Laplacian
of $\bar{\partial_{\cale}}$. We also note that
$\bar{\partial_{\cale}}$ ,$\bar{\partial_{\cale}}^*$ and
$\Delta_{\cale}$ yield
endomorphisms of $\dol{\cale}{ }$. \\

Let  $\Omega_D$ denote the differential form on $C_{k}$ with values
in $\enn(\dol{0}{\cale})$ given by the formula
 $$ \Omega_D = \Pi_{\dol{0}{\cale}} \circ D_k \circ \text{exp}{[-d(t_k - t_{k-1}){\bar{\partial_{\cale}}}^* - (t_k -t_{k-1})\Delta_{\cale}]} \circ ...
\circ  D_1 \circ  {\calg I}_{\dol{0}{\cale}}$$ where
$\Pi_{\dol{0}{\cale}}$ and ${\calg I}_{\dol{0}{\cale}}$ denote the
projection from $\dol{\cale}{ }$ to $\dol{0}{\cale}$ and the
inclusion from $\dol{0}{\cale}$ to $\dol{\cale}{ }$ respectively. \\

As noted in [FLS], to write $\Omega_D$ this way , we require that
$\Delta_{\cale}$ have discrete non negative spectrum (which is the
case for a compact complex manifold). Further, as noted in [FLS],
$\Omega_D$ extends to a $\enn(\dol{0}{\cale})$-valued
(non-homogenous) differential form on $\overline{C_k}$.   \\

We define $${\calg F}_k(D) = \int_{C_{k}} \Omega_D  =
\int_{\overline{C_k}} \Omega_D \text{ . } $$ By the integral over
$\overline{C_k}$ of a non-homogenous differential form, we mean the
integral over $\overline{C_k}$ of its component of top De-Rham
degree. It helps to view the differential forms above as
differential forms on $\overline{C_k}$ rather than on $C_k$ as that
will ensure that the integrals defining the ${\calg F}_k$'s
converge.\\

 Proving that that ${\calg F}_k$ defined in this manner
are the Taylor coefficients of an $\ainf$
morphism is done in [FLS]. \\

 Before we proceed further, we note that if $k=1$, then $C_1$ is a point. The formula for ${\calg F}_k$ given here yields
 $${\calg F}_1(D) = \Pi_{\dol{0}{\cale}} \circ D \circ {\calg I}_{\dol{0}{\cale}} $$ for $D \in \diffdt(\cale)$. In particular,
 if ${\calg D} \in \diff(\cale)$ and if ${\calg D}_* \in \enn( \text{H}^*(X,\cale))$ is the endomorphism of $\text{H}^{*}(X,\cale) $ induced by $\calg D$, then
 ${\calg F}_1(\calg D) = {\calg D}_*$. It follows from Fact 2, Section 3.1 that $\text{tr}([\calg D]) = \text{str}(\calg D_*)$.
 This proves that the following diagram commutes. \\

 $$\begin{CD}
 \hoch{0}(\diff(\cale)) @>>> \hoch{0}(\diffdt(\cale)) \\
 @V{\text{str}}VV   @VV{\text{tr}}V \\
 \compl  @>>{\text{id}}> \compl \\
 \end{CD} $$

 \subsubsection{Rewriting the formula for ${\calg F}_k$.}

 Given an endomorphism $\varphi$ of $\dol{\cale}{ }$, let
${[\varphi]}_i$ denote the endomorphism $id \otimes .... \otimes
\varphi \otimes .. \otimes id$ of $\dol{\cale}{ }^{\otimes k}$ with
$\varphi$ acting on the $i$th
factor from the right. \\

Let $$\Phi = {[id]}_k \circ {[\text{exp}{[-d(t_k -
t_{k-1}){\bar{\partial_{\cale}}}^* - (t_k
-t_{k-1})\Delta_{\cale}]}]}_{k-1} \circ ... $$ $$... \circ
{[\text{exp}{[-d(t_2 - t_{1}){\bar{\partial_{\cale}}}^* - (t_2
-t_{1})\Delta_{\cale}]}]}_1 \text{ . }$$ This is a differential form
on $C_k$ with values in $\enn({\dol{\cale}{ }}^{\otimes k}) $ though
it is not a differential operator. Similarly, if $D= D_1 \otimes ..
\otimes D_k \in \diffdt(\cale)^{\otimes k}$, $D$ yields an
endomorphism $D:= {[D_1]}_k \circ ... \circ {[D_k]}_1$ of
$\enn({\dol{\cale}{ }}^{\otimes k})$.
\\

We have a composition map $m: \enn(\dol{\cale}{ })^{\otimes k} \rar
\enn(\dol{\cale}{ })$. Identifying $\enn({\dol{\cale}{ }}^{\otimes
k})$ with $\enn(\dol{\cale}{ })^{\otimes k}$, we obtain a
composition map $m:\enn({\dol{\cale}{ }}^{\otimes k}) \rar
\enn(\dol{\cale}{ })$. We recall from [FLS] that the formula for $\Omega_D$ can be rewritten as follows \\

$$ \Omega_D = \Pi_{\dol{0}{\cale}} \circ m(\Phi \circ D) \circ {\calg I}_{\dol{0}{\cale}} \text{ . }$$ Thus,
 $${\calg F}_k(D)= \int_{C_k} \Pi_{\dol{0}{\cale}} \circ m(\Phi \circ D) \circ {\calg I}_{\dol{0}{\cale}} \text{ . } $$

\subsection{Extending the supertrace-I.}

Recall that $\dif(\cale)$ denotes the sheaf of holomorphic
differential operators on $\cale$. Let $\dif(\cale)(U)$ denote
$\Gamma(U,\dif(\cale))$ for any open $U \subset X$. Let
$\cale^{\boxtimes k}$ denote the $k$-fold exterior tensor power of
$\cale$ on $X^k:= X \times ... \times X$. We observe that the
differential on the Hochschild complex
$\text{C}^{\bullet}(\dif(\cale)(U))$ extends to a differential on
the graded vector space $\oplus_{k \geq 1} \dif(\cale^{\boxtimes
k}(U^k))[k-1]$. The resulting complex is called the completed
Hochschild
complex of $\dif(\cale)(U)$ and denoted by $\widehat{\hcc{\bullet}{\dif(\cale)(U)}}$. \\

{\it Definition 4:} The completed Hochschild complex of
$\dif(\cale)$ is the sheaf of complexes associated to the preshreaf
$U \leadsto \widehat{\hcc{\bullet}{\dif(\cale)(U)}}$ of complexes of
$\compl$-vector spaces. \\

The completed Hochschild complex of $\dif(\cale)$ is a sheaf of complexes of $\compl$-vector
spaces on $X$. It is denoted in this paper by $\widehat{\text{hoch}}(\dif(\cale))$.\\

{\it Definition 5:} The $i$th completed Hochschild homology
$\choch{i}(\dif(\cale))$ of $\dif(\cale)$ is the hypercohomology
$\hcoh{i}{\widehat{\text{hoch}}(\dif(\cale))}$ of $\widehat{\text{hoch}}(\dif(\cale))$. \\

 Recall that the $\ainf$-map ${\calg F}$ whose construction we recalled in section
3.1 enables us to pull back the supertrace on $\enn(\dol{\cale}{ })$
to a Hochschild $0$-cocycle of $\diffdt(\cale)$. Unfortunately, the
$\ainf$-map $\calg F$ does not automatically enable us to directly
pull back the supertrace to a $\compl$-linear functional on
$\choch{0}(\dif(\cale))$. This subsection is devoted to an important
intermediate step that enables us to construct a $\compl$-linear
functional on $\choch{0}(\dif(\cale))$ which extends
the supertrace  on $\diff(\cale)$. \\

\subsubsection{``Estimating" $\fls$.}

We now use a homotopy very similar to the homotopy in Prop 3.1 [Tel]
to show that the
FLS functional of any $0$-cycle in $\text{C}^{\bullet}(\diffdt(\cale))$ depends only on its component in $\text{Diff}^{0}(\cale)$.\\

\textbf{Construction 1:}Let $\phi:X \times X \rar [0,\infty)$ be a
Riemannian distance. Let $t$ be any positive real number. We can
choose a finite cover of $X$ by open sets $U_i , 1 \leq i \leq m$
such that $\phi(x,y) < t$ for all $x,y \in U_i$ for any $i$. Choose
a partition of unity $\{f_i\}$ by compactly supported (nonnegative
valued) smooth functions subordinate to the cover $\{U_i\}$. Let
$g_i$ be a compactly supported smooth function on $U_i$ with values
in $[0,1]$ that is identically $1$ on the support of $f_i$. Then,
$f:= \sum_{i=1}^{i=m}f_i \otimes g_i$ is a smooth function on $X
\times X$ whose restriction to the diagonal is identically $1$.
Also, $f$ vanishes outside the subset $\{\phi(x,y) \leq t\}$ of $X
\times X$, and the maximum value of $f$ on $X \times
X$ is $1$. \\

Somewhat as in [Tel], let $E_{f}:\diffdt(\cale)^{\otimes k} \rar
\diffdt(\cale)^{\otimes k+1}$ be the map
$$D_1 \otimes ... \otimes D_k \mapsto -\sum_{i=1}^{i=m} f_i \otimes g_iD_1 \otimes D_2 \otimes .... \otimes D_k \text{ . }$$
Let $d$ be the differential of the complex
$\text{C}^{\bullet}(\diffdt(\cale))$. Then,

$$dE_f + E_fd = 1 -N_f $$

where
$$
N_f(D_1 \otimes ... \otimes D_k) = \pm \sum_{i=1}^{i=m} \bar\partial
f_i \otimes g_iD_1 \otimes ... \otimes D_k \pm \sum_{i=1}^{i=m} f_i
\otimes (\bar{\partial}g_i)D_1 \otimes ... \otimes D_k $$

\begin{equation} \label{main} \pm \sum_{i=1}^{i=m} f_i \otimes
g_iD_kD_1 \otimes .... \otimes D_{k-1} \pm \sum_{i=1}^{i=m} D_{k}f_i
\otimes g_iD_1 \otimes ... \otimes D_{k-1} \text{    } \forall
\text{   } k \geq 2
\text{ . }\\
\end{equation}

\textbf{Basic argument.} Note that if $\alpha$ is a cycle in
$\text{C}^{\bullet}(\diffdt(\cale))$, then $\alpha$ is homotopic to
$N_f\alpha$. This is true for all $t >0$. Also, the FLS linear
functional $I_{\text{FLS}}$ is a Hochschild $0$-cocycle of
$\diffdt(\cale)$. It follows that
$$I_{\text{FLS}}(\alpha) = I_{\text{FLS}}(N_f\alpha) \text{ . }$$
Let $\alpha=\alpha_1+...+\alpha_p$ with $\alpha_i \in
\diffdt(\cale)^{\otimes i}$. We then show that
$$\sum_{i \geq 2} |\fls(N_f\alpha_i)| \leq C\epsilon(t)
$$ where $C$ is a constant that only depends on $\alpha$ and
$\epsilon(t).\text{Vol}(X \times X)$ is the volume of the subset
$\{\phi(x,y) \leq t\}$ of $X \times X$. Similarly, we show that
$|\fls(\alpha_1)-\fls(N_f(\alpha_1))| \leq C'\epsilon(t)$ where $C'$

depends only on $\alpha$. Since the construction of $f$ as in
Construction 1 is possible for all positive $t$, letting $t$
approach $0$ we see that
$I_{\text{FLS}}(\alpha)= \fls(\alpha_1)$. \\

\textbf{Remark:} Recall that $\dol{0}{\cale}$ can be identified with
the kernel of the Laplacian $\Delta_{\cale}$. Let $\dol{L^2}{\cale}$
denote the Hilbert space of square integrable sections of
$\text{Dolb}^{\bullet}(X,\strc) \otimes_{\strc} \dif(\cale)$. Then
${\calg I}_{\dol{0}{\cale}} \circ \Pi_{\dol{0}{\cale}}$ is an
integral operator on $\dol{L^2}{\cale}$ with smooth kernel that
projects onto the image of ${\calg I}_{\dol{0}{\cale}}$ (see [BGV]
Chapter 2). We will denote this operator by $\Pi_{0,\cale}$ or
$\Pi_0$ when there is no confusion regarding the vector bundle being
used. Also, one can check that if $M \in \text{End}(\dol{}{\cale})$,
then $\Pi_0 M\Pi_0$ makes sense as a trace class operator on
$\dol{L^2}{\cale}$, and
$$\text{str}(\Pi_0M\Pi_0) = \text{str}(\Pi_{\dol{0}{\cale}} \circ M
\circ {\calg I}_{\dol{0}{\cale}}) \text{ . }$$

\begin{prop}Let $D_1,...,D_k \in \diffdt(\cale)$. Then,
$$|\fls(N_f(D_1 \otimes ... \otimes D_k))| \leq C\epsilon(t)
\text{ . }$$ The constant $C$ above depends only on $D_1,..,D_k$.\\
\end{prop}

\begin{proof}

{\it Part 1: The setup.}\\

Recall that

$$
N_f(D_1 \otimes ... \otimes D_k) = \pm \sum_{i=1}^{i=m} \bar\partial
f_i \otimes g_iD_1 \otimes ... \otimes D_k \pm \sum_{i=1}^{i=m} f_i
\otimes (\bar{\partial}g_i)D_1 \otimes ... \otimes D_k $$

$$ \pm \sum_{i=1}^{i=m} f_i \otimes g_iD_kD_1 \otimes .... \otimes
D_{k-1} \pm \sum_{i=1}^{i=m} D_{k}f_i \otimes g_iD_1 \otimes ...
\otimes D_{k-1} \text{    } \forall \text{   } k \geq 2
\text{ . }\\
$$

We estimate the FLS functional of each of the summands on the RHS
separately. Let ${\calg F}_k$ denote the $k$-th Taylor component of
the FLS $A_{\infty}$ map. Let ${C_k}$ denote the configuration space
$\Pi_{i=1}^{k-1} \{\tau_i > 0\}$. Then, if $\alpha_i \in
\diffdt(\cale)$,
$${\calg F}_k(\alpha_1 \otimes ... \otimes \alpha_k) = \int_{C_k}
\Pi_0 \alpha_1 \bar{\partial}^* \text{e}^{-\tau_1 \Delta} \alpha_2
.. \text{e}^{-\tau_{k-1}\Delta} \alpha_k \Pi_0 d\tau_1...d\tau_{k-1}
\text{ . } $$

Let $Q_{\tau}(x,y) \in \Gamma(X \times X, \cale \otimes
\Omega^{0,\bullet} \boxtimes {(\Omega^{0,\bullet})}^* \otimes
\cale^*) $ denote the (smooth) kernel of the operator
$\bar{\partial}^*\text{e}^{-\tau\Delta}$. Let $p_{\infty}(x,y)$
denote the (smooth) kernel of $\Pi_0$. Then, the kernel of ${\calg
F}_k(\alpha_1 \otimes ... \otimes \alpha_k)$ is
$$\int_{C_k} \int_{X^k}
p_{\infty}(y_0,x_1)\alpha_{1,x_1}Q_{\tau_1}(x_1,x_2)....
Q_{\tau_{k-1}}\alpha_{k,x_k}p_{\infty}(x_k,y_0') |dx_1|...|dx_k|
d\tau_1...d\tau_{k-1} \text{ . }$$

Recall that $$\fls(\alpha_1 \otimes ... \otimes \alpha_k) =
\sum_{s=0}^{s=k-1} \text{str}{\calg F}_k(((\sigma^s(\alpha))))
\text{ . }$$ Here $\sigma$ is a (signed) cyclic permutation. Each of
the $k$ cyclic permutations of $\alpha_1 \otimes .. \otimes
\alpha_k$ yields exactly one summand contributing towards
$\fls(\alpha_1
\otimes... \otimes \alpha_k)$. \\

{\it Part 2: Estimating $\fls(\pm \sum_{i=1}^{i=m} \bar\partial
f_i \otimes g_iD_1 \otimes ... \otimes D_k) $ .}\\

The cyclic permutations of $\sum_{i=1}^{i=m} \bar\partial
f_i \otimes g_iD_1 \otimes ... \otimes D_k $ are as follows :\\
$$\pm \sum_{i=1}^{i=m} \bar\partial f_i \otimes g_iD_1 \otimes ... \otimes D_k $$
$$\pm \sum_{l=1}^{l=m} D_{k-i+1} \otimes ... \otimes D_k \otimes  \bar\partial f_l \otimes g_lD_1 \otimes ...\otimes D_{k-i} \text{  } 1 \leq i \leq k-1 $$
$$ \pm \sum_{l=1}^{l=m} g_lD_1 \otimes .... \otimes D_k \otimes  \bar\partial f_l \text{ . }$$

We shall henceforth denote the kernel of an integral operator $T$ on
$\Gamma(X,\cale \otimes \Omega^{0,\bullet}(X))$ by ${\overline{T}}$.
Let $D^a$ denote the formal Adjoint of a differential operator $D$
on $\Gamma(X,\cale \otimes \Omega^{0,\bullet}(X))$. This is a
differential operator on $\Gamma(X,\Omega^{0,\bullet}(X)^* \otimes
\cale^*)$ since the bundle of densities on $X$ has a canonical
trivialization (see [BGV] Chapter 2). With this in mind,

$$\overline{{\calg F}_{k+1}}( \sum_{i=1}^{i=m} \bar\partial
f_i \otimes g_iD_1 \otimes ... \otimes D_k) $$
$$ = \int_{C_{k+1}} \int_{X^{k+1}} \sum_{i=1}^{i=m}
p_{\infty}(y_0,x_0) (\bar{\partial} f_i)(x_0)
Q_{\tau_1}(x_0,x_1)g_i(x_1)D_{1,x_1}Q_{\tau_1}(x_1,x_2)...$$
$$ ....p_{\infty}(x_k,y_0') |dx_0|...|dx_k| d\tau_1...d\tau_k $$
$$ = \pm \int_{C_{k+1}} \int_{X^{k+1}} \sum_{i=1}^{i=m}
 p_{\infty}(y_0,x_0) \bar{\partial}_{x_0} f_i(x_0)
Q_{\tau_1}(x_0,x_1) g_i(x_1)D_{1,x_1}Q_{\tau_1}(x_1,x_2)...$$
$$ ....p_{\infty}(x_k,y_0') |dx_0|...|dx_k| d\tau_1...d\tau_k $$
$$ \pm \int_{C_{k+1}} \int_{X^{k+1}} \sum_{i=1}^{i=m}
p_{\infty}(y_0,x_0) f_i(x_0) \bar{\partial}_{x_0}
Q_{\tau_1}(x_0,x_1) g_i(x_1)D_{1,x_1}Q_{\tau_1}(x_1,x_2)...$$
$$ ....p_{\infty}(x_k,y_0') |dx_0|...|dx_k| d\tau_1...d\tau_k $$
$$ = \pm \int_{C_{k+1}} \int_{X^{k+1}}
[\bar{\partial}^a_{x_0}p_{\infty}(y_0,x_0)] [ \sum_{i=1}^{i=m}
f_i(x_0)g_i(x_1)] [Q_{\tau_1}(x_0,x_1)][
D_{1,x_1}Q_{\tau_1}(x_1,x_2)]...$$
$$ ....[...p_{\infty}(x_k,y_0')] |dx_0|...|dx_k| d\tau_1...d\tau_k $$
$$ \pm \int_{C_{k+1}} \int_{X^{k+1}}
[p_{\infty}(y_0,x_0)] [\sum_{i=1}^{i=m}
f_i(x_0)g_i(x_1)][\bar{\partial}_{x_0}
Q_{\tau_1}(x_0,x_1)][D_{1,x_1}Q_{\tau_1}(x_1,x_2)]...$$
\begin{equation} \label{case1}....[...p_{\infty}(x_k,y_0')] |dx_0|...|dx_k| d\tau_1...d\tau_k  \end{equation}

The last equality in $\eqref{case1}$ is valid since $X$ is compact.
We also use the fact that if $h(x,y) \in \Gamma(X \times X , \cale
\otimes \Omega^{0,\bullet} \boxtimes {\Omega^{0,\bullet}}^* \otimes
\cale^*)$, then $a(x)h(x,y)b(y) = \pm a(x)b(y)h(x,y)$ for smooth
functions $a,b$ on $X$. We now note that if $p_{\tau}$ denotes the
kernel of $e^{-\tau\Delta}$ then $$Q_{\tau}(x,y) =
\bar{\partial}^*_xp_{\tau}(x,y) = \bar{\partial}^*
p_{(0,\infty)}(x)p_{\tau}(x,y)p_{(0,\infty)}(y) $$ where
$p_{(0,\infty)}$ is the projection to the span of the eigenvectors
of $\Delta$ corresponding to positive eigenvalues as in [BGV]
Proposition 2.37. This is because ${\bar{\partial}}^*$ kills the
kernel of $\Delta$. It follows from [BGV] Proposition 2.37 that
\begin{equation} \label{estt} ||Q_{\tau}(x,y)||_l \leq
C(||.||_l)\text{e}^{-\frac{1}{2}\lambda_1\tau}\end{equation}  where
$||.||_l$ is any $C^l$ norm on $\Gamma(X \times X,\cale \otimes
\Omega^{0,\bullet} \boxtimes {\Omega^{0,\bullet}}^* \otimes
\cale^*)$ and
$\lambda_1$ is the first positive eigenvalue of $\Delta_{\cale}$. \\

It follows from $\eqref{estt}$ that the sup-norm of each term within
a square bracket in each integrand that involves a $Q_{\tau_i}$ is
bounded above by $C_i\text{e}^{-\tau_i\frac{\lambda_1}{2}}$ for some
constant $C_i$ depending only on $D_1 \otimes ... \otimes D_k$. The
remaining terms have finite sup-norm. The sup-norm of
$\sum_{i=1}^{i=m} f_i(x_0)g_i(x_1)$ is $1$. Hence, there exists a
positive constant $C$ depending only on $D_1 \otimes ... \otimes
D_k$ such that $$ ||\overline{{\calg F}_{k+1}}(\sum_{i=1}^{i=m}
\bar\partial f_i \otimes g_iD_1 \otimes ... \otimes D_k)||_0 \leq
C\int_{C_{k+1}} \text{e}^{-\tau_1\frac{\lambda_1}{2}}.$$
\begin{equation} \label{estt1}...\text{e}^{-\tau_k\frac{\lambda_k}{2}} \epsilon(t) \text{Vol}(X)^{k+1}
d\tau_1...d\tau_k  =  \frac{2^k}{\lambda_1^k} C\epsilon(t)
\text{Vol}(X)^{k+1} \text{ . }\end{equation} where
$\epsilon(t)\text{Vol}(X \times X)$ is the volume of the support of
$\{\phi(x,y) \leq t\}$ in $X \times X$. Now, $$ \text{str}({\calg
F}_{k+1}(\pm \sum_{i=1}^{i=m} \bar\partial f_i \otimes g_iD_1
\otimes ... \otimes D_k))$$ $$ = \pm \int_X
\text{str}(\overline{\calg F_{k+1}}(\sum_{i=1}^{i=m} \bar\partial
f_i \otimes g_iD_1 \otimes ... \otimes D_k)(y_0,y_0)|dy_0|$$ It
follows from $\eqref{estt1}$ that
$$ |\text{str}({\calg F}_{k+1}(\pm \sum_{i=1}^{i=m}
\bar\partial f_i \otimes g_iD_1 \otimes ... \otimes D_k))| \leq
C_1\epsilon(t) $$ for some constant $C_1$ depending on $D_1 \otimes
... \otimes D_k$ only.\\

The same method is used to estimate the other summands that
contribute to $\fls(\pm \sum_{i=1}^{i=m} \bar\partial f_i \otimes
g_iD_1 \otimes ... \otimes D_k))$ and show that those contributions
are at most $C'\epsilon(t)$ as well - in fact the FLS functionals of
the other summands of $N_f(D_1 \otimes ... \otimes D_k)$ are
estimated in the same way. We however, present the detailed
calculations in the next 2 parts to be extra careful.\\

{\it Part 3: Estimating $\fls(\pm \sum_{i=1}^{i=m} \bar\partial
f_i \otimes g_iD_1 \otimes ... \otimes D_k))$ -II.}\\

We show the equivalent of the calculation $\eqref{case1}$ for the
other summands contributing to $\fls(\pm \sum_{i=1}^{i=m}
\bar\partial f_i \otimes g_iD_1 \otimes ... \otimes D_k)$. At the
end of each calculation, the integrand in the integral computing the
kernel of the operator whose supertrace we need is written as a
product of terms marked by square brackets. The argument in Part 2
then works almost word for word to show that the corresponding
contribution to $\fls(N_f(D_1 \otimes ... \otimes D_k))$ is bounded
by a constant depending on the $D_i$'s times $\epsilon(t)$.  \\

$$\overline{{\calg F}_{k+1}}( \sum_{l=1}^{l=m} D_{k-i+1} \otimes ... \otimes D_k \otimes \bar{\partial} f_l \otimes g_lD_1 \otimes ...\otimes D_{k-i}) $$
$$ = \int_{C_{k+1}} \int_{X^{k+1}} \sum_{l=1}^{l=m}
p_{\infty}(y_0,x_0)D_{k-i+1,x_0}Q_{\tau_1}(x_0,x_1).....$$
$$ .....D_{k,x_{i-1}}Q_{\tau_i}(x_{i-1},x_i)(\bar\partial
f_l)(x_i)Q_{\tau_{i+1}}(x_i,x_{i+1})g_l(x_{i+1})D_{1,x_{i+1}}.....$$
$$.....p_{\infty}(x_k,y_0') |dx_0|...|dx_k| d\tau_1...d\tau_k $$
$$ = \pm \int_{C_{k+1}} \int_{X^{k+1}} \sum_{l=1}^{l=m}
p_{\infty}(y_0,x_0)D_{k-i+1,x_0}Q_{\tau_1}(x_0,x_1).....$$
$$ .....D_{k,x_{i-1}}Q_{\tau_i}(x_{i-1},x_i)\bar\partial_{x_i}
f_l(x_i)Q_{\tau_{i+1}}(x_i,x_{i+1})g_lx_{i+1})D_{1,x_{i+1}}.....$$
$$.....p_{\infty}(x_k,y_0') |dx_0|...|dx_k| d\tau_1...d\tau_k $$
$$ \pm \int_{C_{k+1}} \int_{X^{k+1}} \sum_{l=1}^{l=m}
p_{\infty}(y_0,x_0)D_{k-i+1,x_0}Q_{\tau_1}(x_0,x_1).....$$
$$ .....D_{k,x_{i-1}}Q_{\tau_i}(x_{i-1},x_i)f_l(x_i)\bar\partial_{x_i}Q_{\tau_{i+1}}(x_i,x_{i+1})
g_l(x_{i+1})D_{1,x_{i+1}}.....$$
$$.....p_{\infty}(x_k,y_0') |dx_0|...|dx_k| d\tau_1...d\tau_k $$
$$ = \pm \int_{C_{k+1}} \int_{X^{k+1}}
[p_{\infty}(y_0,x_0)][D_{k-i+1,x_0}Q_{\tau_1}(x_0,x_1)].....$$
$$ .....[D_{k,x_{i-1}}\bar\partial^a_{x_i}Q_{\tau_i}(x_{i-1},x_i)]
[\sum_{l=1}^{l=m}f_l(x_i)g_l(x_{i+1})][Q_{\tau_{i+1}}(x_i,x_{i+1})][D_{1,x_{i+1}}..]...$$
$$...[..p_{\infty}(x_k,y_0')] |dx_0|...|dx_k| d\tau_1...d\tau_k $$
$$ \pm \int_{C_{k+1}} \int_{X^{k+1}}
[p_{\infty}(y_0,x_0)][D_{k-i+1,x_0}Q_{\tau_1}(x_0,x_1)].....$$
$$
.....[D_{k,x_{i-1}}Q_{\tau_i}(x_{i-1},x_i)][\sum_{l=1}^{l=m}
f_l(x_i)g_l(x_{i+1})][\bar\partial_{x_i}Q_{\tau_{i+1}}(x_i,x_{i+1})]
[D_{1,x_{i+1}}..]...$$
\begin{equation} \label{case2} ...[...p_{\infty}(x_k,y_0')] |dx_0|...|dx_k| d\tau_1...d\tau_k
\end{equation}

$${\calg F}_{k+1}( \sum_{l=1}^{l=m}g_lD_1 \otimes .... \otimes D_k \otimes \bar\partial
f_l) $$
$$ = \int_{C_{k+1}} \int_{X^{k+1}}\sum_{l=1}^{l=m} p_{\infty}(y_0,x_0)g_l(x_0)D_{1,x_0}Q_{\tau_1}(x_0,x_1).....
$$ $$....D_kQ_{\tau_k}(x_{k-1},x_k)(\bar\partial
f_l)(x_k)p_{\infty}(x_k,y_0') |dx_0|...|dx_k|d\tau_1...d\tau_k$$
$$ = \pm \int_{C_{k+1}} \int_{X^{k+1}} \sum_{l=1}^{l=m}
p_{\infty}(y_0,x_0) g_l(x_0)D_{1,x_0}Q_{\tau_1}(x_0,x_1).....
$$ $$....D_kQ_{\tau_k}(x_{k-1},x_k)\bar\partial_{x_k}
f_l(x_k)p_{\infty}(x_k,y_0') |dx_0|...|dx_k|d\tau_1...d\tau_k$$
$$  \pm \int_{C_{k+1}} \int_{X^{k+1}} \sum_{l=1}^{l=m}
p_{\infty}(y_0,x_0) g_l(x_0)D_{1,x_0}Q_{\tau_1}(x_0,x_1).....
$$ $$....D_kQ_{\tau_k}(x_{k-1},x_k)
f_l(x_k)\bar\partial_{x_k}p_{\infty}(x_k,y_0')
|dx_0|...|dx_k|d\tau_1...d\tau_k$$
$$ =  \pm \int_{C_{k+1}} \int_{X^{k+1}}
[p_{\infty}(y_0,x_0)] [\sum_{l=1}^{l=m}
f_l(x_k)g_l(x_0)][D_{1,x_0}Q_{\tau_1}(x_0,x_1)].....
$$ $$....[D_k\bar\partial^a_{x_k}Q_{\tau_k}(x_{k-1},x_k)]
[p_{\infty}(x_k,y_0')] |dx_0|...|dx_k|d\tau_1...d\tau_k$$
$$  \pm \int_{C_{k+1}} \int_{X^{k+1}}
[p_{\infty}(y_0,x_0)] [\sum_{l=1}^{l=m}
f_l(x_k)g_l(x_0)][D_{1,x_0}Q_{\tau_1}(x_0,x_1)].....
$$ \begin{equation} \label{case3}...[D_kQ_{\tau_k}(x_{k-1},x_k)]
[\bar\partial_{x_k}p_{\infty}(x_k,y_0')]
|dx_0|...|dx_k|d\tau_1...d\tau_k \end{equation}

{\it Part 4: Estimating other summands contributing to $\fls(N_f(D_1
\otimes .... \otimes D_k))$.}\\

Recall that
$$
N_f(D_1 \otimes ... \otimes D_k) = \pm \sum_{i=1}^{i=m} \bar\partial
f_i \otimes g_iD_1 \otimes ... \otimes D_k \pm \sum_{i=1}^{i=m} f_i
\otimes (\bar{\partial}g_i)D_1 \otimes ... \otimes D_k $$

$$ \pm \sum_{i=1}^{i=m} f_i \otimes g_iD_kD_1 \otimes .... \otimes
D_{k-1} \pm \sum_{i=1}^{i=m} D_{k}f_i \otimes g_iD_1 \otimes ...
\otimes D_{k-1} \text{    } \forall \text{   } k \geq 2
\text{ . }\\
$$

The cyclic permutations of $\sum_{l=1}^{l=m}f_l \otimes g_lD_kD_1
\otimes .... \otimes D_{k-1}$ that contribute to its FLS functional
are as follows:
$$\sum_{l=1}^{l=m}f_l \otimes g_lD_kD_1 \otimes .... \otimes D_{k-1}$$
$$D_{k-i} \otimes .... \otimes D_{k-1} \otimes
\sum_{l=1}^{l=m}f_l \otimes g_lD_kD_1 \otimes ... \otimes D_{k-i-1}
\text{  for }1\leq i \leq k-2$$
$$ \sum_{l=1}^{l=m}g_lD_kD_1 \otimes....\otimes D_{k-1} \otimes f_l $$

As in earlier parts,
$$\overline{{\calg F}_k}(\sum_{l=1}^{l=m}f_l \otimes g_lD_kD_1 \otimes .... \otimes
D_{k-1})$$
$$= \int_{C_k} \int_{X^k}
[p_{\infty}(y_0,x_1)][\sum_{l=1}^{l=m}f_l(x_1)
g_l(x_2)][Q_{\tau_1}(x_1,x_2)][D_{k,x_2}D_{1,x_2}...]$$
\begin{equation} \label{case4}...[D_{k-1}p_{\infty}(x_k,y_0')]
|dx_1|...|dx_k| d\tau_1...d\tau_{k-1} \text{ . } \end{equation}
$$\overline{{\calg F}_k}(D_{k-i} \otimes .... \otimes D_{k-1} \otimes
\sum_{l=1}^{l=m}f_l \otimes g_lD_kD_1 \otimes ... \otimes
D_{k-i-1})$$
$$ = \int_{C_k} \int_{X^k}
[p_{\infty}(y_0,x_1)][D_{k-i,x_1}Q_{\tau_1}(x_1,x_2)]....$$
$$ ...[\sum_{l=1}^{l=m}f_l(x_{i+1})g_l(x_{i+2})]
[Q_{\tau_{i+1}}(x_{i+1},x_{i+2})][D_{k,x_{i+2}}D_{1,x_{i+2}}...]...
$$ \begin{equation} \label{case5} [..p_{\infty}(x_k,y_0')]
 |dx_1|...|dx_k|
d\tau_1...d\tau_{k-1}  \text{ . }\end{equation}
$$\overline{{\calg F}_k}(\sum_{l=1}^{l=m}g_lD_kD_1 \otimes....\otimes D_{k-1} \otimes f_l)$$
$$ =\int_{C_k} \int_{X^k}
[p_{\infty}(y_0,x_1)][\sum_{l=1}^{l=m}f_l(x_k)g_l(x_1)][D_{k,x_1}D_{1,x_1}Q_{\tau_1}(x_1,x_2)].....
$$ \begin{equation} \label{case6} ......[p_{\infty}(x_k,y_0')] |dx_1|...|dx_k| d\tau_1...d\tau_{k-1}
\end{equation}

For each of the summands dealt with by equations $\eqref{case4}$,
$\eqref{case5}$ and $\eqref{case6}$ respectively, one can argue
exactly as we did after $\eqref{case1}$.\\

The cyclic permutations of $\sum_{l=1}^{l=m}D_{k}f_l \otimes g_lD_1
\otimes ... \otimes D_{k-1}$ that contribute to its FLS functional
are
$$\sum_{l=1}^{l=m} D_{k}f_l \otimes g_lD_1 \otimes ... \otimes D_{k-1}$$
$$\sum_{l=1}^{l=m} D_{k-i} \otimes .. \otimes D_{k-1} \otimes
D_{k}f_l \otimes g_lD_{1} \otimes .. \otimes D_{k-i-1} \text{ for }
1 \leq i \leq k-2 $$
$$ \sum_{l=1}^{l=m} g_lD_1 \otimes ... \otimes D_{k}
f_l$$

$$\overline{{\calg F}_k}(D_{k}f_l \otimes g_lD_1 \otimes ... \otimes
D_{k-1})$$
$$ = \int_{C_k} \int_{X^k} \sum_{l=1}^{l=m} p_{\infty}(y_0,x_1)
D_{k,x_1}f_l(x_1)Q_{\tau_1}(x_1,x_2)g_l(x_2)D_{1,x_2}....
$$ $$ ...D_{k-1,x_k}p_{\infty}(x_k,y_0') |dx_1|...|dx_k|
d\tau_1...d\tau_{k-1} $$
$$ = \int_{C_k} \int_{X^k} [D^a_{k,x_1}p_{\infty}(y_0,x_1)]
[\sum_{l=1}^{l=m}
f_l(x_1)g_l(x_2)][Q_{\tau_1}(x_1,x_2)][D_{1,x_2}..]..
$$ \begin{equation} \label{case7} ...[D_{k-1,x_k}p_{\infty}(x_k,y_0')]
|dx_1|...|dx_k| d\tau_1...d\tau_{k-1} \end{equation}

$$\overline{{\calg F}_k}(\sum_{l=1}^{l=m} D_{k-i} \otimes .. \otimes D_{k-1} \otimes
D_{k}f_l \otimes g_lD_{1} \otimes .. \otimes D_{k-i-1})$$
$$ = \int_{C_k} \int_{X^k} \sum_{l=1}^{l=m}
p_{\infty}(y_0,x_1)D_{k-i,x_1}Q_{\tau_1}(x_1,x_2)....$$
$$...D_{k-1,x_i}Q_{\tau_i}(x_i,x_{i+1}) D_{k,x_{i+1}}f_l(x_{i+1})
Q_{\tau_{i+1}}(x_{i+1},x_{i+2})g_l(x_{i+2})D_{1,x_{i+2}}....$$
$$....p_{\infty}(x_k,y_0') |dx_1|...|dx_k|
d\tau_1...d\tau_{k-1} $$
$$ = \int_{C_k} \int_{X^k}
[p_{\infty}(y_0,x_1)][D_{k-i,x_1}Q_{\tau_1}(x_1,x_2)]....$$
$$...[D_{k-1,x_i}D^a_{k,x_{i+1}}Q_{\tau_i}(x_i,x_{i+1})][
\sum_{l=1}^{l=m}f_l(x_{i+1})g_l(x_{i+2})]
[Q_{\tau_{i+1}}(x_{i+1},x_{i+2})][D_{1,x_{i+2}}..]..$$
\begin{equation}\label{case8} ....[..p_{\infty}(x_k,y_0')] |dx_1|...|dx_k| d\tau_1...d\tau_{k-1}
\end{equation}

$$\overline{{\calg F}_k}( \sum_{l=1}^{l=m} g_lD_1 \otimes ... \otimes D_{k}
f_l)$$
$$ = \int_{C_k} \int_{X^k} \sum_{l=1}^{l=m}
p_{\infty}(y_0,x_1)g_l(x_1)D_{1,x_1}Q_{\tau_1}(x_1,x_2)...
$$ $$...D_{k-1,x_{k-1}}Q_{\tau_{k-1}}(x_{k-1},x_k)
D_{k,x_k} f_l(x_k) p_{\infty}(x_k,y_0') |dx_1|...|dx_k|
d\tau_1...d\tau_{k-1}$$
$$ = \int_{C_k} \int_{X^k}
[p_{\infty}(y_0,x_1)][\sum_{l=1}^{l=m}
f_l(x_k)g_l(x_1)][D_{1,x_1}Q_{\tau_1}(x_1,x_2)]...
$$ \begin{equation} \label{case9}...[D_{k-1,x_{k-1}}D^a_{k,x_k}Q_{\tau_{k-1}}(x_{k-1},x_k)]
[p_{\infty}(x_k,y_0')] |dx_1|...|dx_k| d\tau_1...d\tau_{k-1}
\end{equation}

For each of the summands dealt with by equations $\eqref{case7}$,
$\eqref{case8}$ and $\eqref{case9}$ respectively, one can argue
exactly as we did after $\eqref{case1}$. \\

The cyclic permutations of $\sum_{l=1}^{l=m} f_l \otimes
(\bar{\partial}g_l)D_1 \otimes ... \otimes D_k$ that contribute to
its FLS functional are as follows:
$$\sum_{l=1}^{l=m} f_l \otimes
(\bar{\partial}g_l)D_1 \otimes ... \otimes D_k$$
$$ \sum_{l=1}^{l=m} D_{k-i} \otimes .. \otimes D_k \otimes f_l \otimes
(\bar{\partial}g_l)D_1 \otimes .. \otimes D_{k-i-1} \text{ for } 0
\leq i \leq k-2 $$
$$\sum_{l=1}^{l=m} (\bar{\partial}g_l)D_1 \otimes ... \otimes D_k
\otimes f_l \text{ . }$$

$$\overline{{\calg F}_{k+1}}(\sum_{l=1}^{l=m} f_l \otimes
(\bar{\partial}g_l)D_1 \otimes ... \otimes D_k) $$
$$ \int_{C_{k+1}} \int_{X^{k+1}} \sum_{l=1}^{l=m}
p_{\infty}(y_0,x_0)
f_l(x_0)Q_{\tau_1}(x_0,x_1)(\bar{\partial}g_l)(x_1)D_{1,x_1}...$$
$$...p_{\infty}(x_k,y_0') |dx_0|....|dx_k|d\tau_1...d\tau_k $$
$$= \pm \int_{C_{k+1}} \int_{X^{k+1}} \sum_{l=1}^{l=m}
p_{\infty}(y_0,x_0)
f_l(x_0)Q_{\tau_1}(x_0,x_1)\bar{\partial}_{x_1}g_l(x_1)D_{1,x_1}...$$
$$...p_{\infty}(x_k,y_0') |dx_0|....|dx_k|d\tau_1...d\tau_k $$
$$ \pm \int_{C_{k+1}} \int_{X^{k+1}} \sum_{l=1}^{l=m}
p_{\infty}(y_0,x_0)
f_l(x_0)Q_{\tau_1}(x_0,x_1)g_l(x_1)\bar{\partial}_{x_1}D_{1,x_1}...$$
$$...p_{\infty}(x_k,y_0') |dx_0|....|dx_k|d\tau_1...d\tau_k $$
$$= \pm \int_{C_{k+1}} \int_{X^{k+1}}
[p_{\infty}(y_0,x_0)]
[\sum_{l=1}^{l=m}f_l(x_0)g_l(x_1)][\bar{\partial}^a_{x_1}Q_{\tau_1}(x_0,x_1)][D_{1,x_1}..].$$
$$..[..p_{\infty}(x_k,y_0')] |dx_0|....|dx_k|d\tau_1...d\tau_k$$
$$ \pm \int_{C_{k+1}} \int_{X^{k+1}}
[p_{\infty}(y_0,x_0)] [\sum_{l=1}^{l=m}
f_l(x_0)g_l(x_1)][Q_{\tau_1}(x_0,x_1)][\bar{\partial}_{x_1}D_{1,x_1}..]...$$
\begin{equation} \label{case10}..[..p_{\infty}(x_k,y_0')]|dx_0|....|dx_k|d\tau_1...d\tau_k
\end{equation}

$$\overline{{\calg F}_{k+1}}(\sum_{l=1}^{l=m} D_{k-i} \otimes .. \otimes D_k \otimes f_l \otimes
(\bar{\partial}g_l)D_1 \otimes .. \otimes D_{k-i-1})$$
$$= \int_{C_{k+1}} \int_{X^{k+1}} \sum_{l=1}^{l=m}
p_{\infty}(y_0,x_0)......Q_{\tau_i}(x_{i-1},x_i)f_l(x_i)Q_{\tau_{i+1}}(x_i,x_{i+1})(\bar{\partial}g_l)(x_{i+1})D_{1,x_{i+1}}..$$
$$...p_{\infty}(x_k,y_0') |dx_0|....|dx_k|d\tau_1...d\tau_k $$
$$= \pm \int_{C_{k+1}} \int_{X^{k+1}} \sum_{l=1}^{l=m}
p_{\infty}(y_0,x_0)......Q_{\tau_i}(x_{i-1},x_i)f_l(x_i)Q_{\tau_{i+1}}(x_i,x_{i+1})\bar{\partial}_{x_{i+1}}g_l(x_{i+1})D_{1,x_{i+1}}..$$
$$...p_{\infty}(x_k,y_0') |dx_0|....|dx_k|d\tau_1...d\tau_k $$
$$ \pm \int_{C_{k+1}} \int_{X^{k+1}} \sum_{l=1}^{l=m}
p_{\infty}(y_0,x_0)......Q_{\tau_i}(x_{i-1},x_i)f_l(x_i)Q_{\tau_{i+1}}(x_i,x_{i+1})g_l(x_{i+1})\bar{\partial}_{x_{i+1}}D_{1,x_{i+1}}..$$
$$...p_{\infty}(x_k,y_0') |dx_0|....|dx_k|d\tau_1...d\tau_k $$
$$= \pm \int_{C_{k+1}} \int_{X^{k+1}}
[p_{\infty}(y_0,x_0)]....[..Q_{\tau_i}(x_{i-1},x_i)][\sum_{l=1}^{l=m}f_l(x_i)g_l(x_{i+1})][\bar{\partial}^a_{x_{i+1}}Q_{\tau_{i+1}}(x_i,x_{i+1})]$$
$$ [D_{1,x_{i+1}}...]...[...p_{\infty}(x_k,y_0')] |dx_0|....|dx_k|d\tau_1...d\tau_k $$
$$= \pm \int_{C_{k+1}} \int_{X^{k+1}}
[p_{\infty}(y_0,x_0)]....[..Q_{\tau_i}(x_{i-1},x_i)][\sum_{l=1}^{l=m}f_l(x_i)g_l(x_{i+1})][Q_{\tau_{i+1}}(x_i,x_{i+1})]$$
\begin{equation} \label{case11} [\bar{\partial}_{x_{i+1}}D_{1,x_{i+1}}..
.]...[...p_{\infty}(x_k,y_0')] |dx_0|....|dx_k|d\tau_1...d\tau_k
\end{equation}

$$\overline{{\calg F}_{k+1}}(\sum_{l=1}^{l=m} (\bar{\partial}g_l)D_1 \otimes ... \otimes D_k
\otimes f_l ) $$
$$ = \int_{C_{k+1}} \int_{X^{k+1}} \sum_{l=1}^{l=m}
p_{\infty}(y_0,x_0)(\bar{\partial}g_l)(x_0)D_{1,x_0}Q_{\tau_1}(x_0,x_1)......$$
$$....f_l(x_k)p_{\infty}(x_k,y_0') |dx_0|....|dx_k|d\tau_1...d\tau_k $$
$$ = \pm \int_{C_{k+1}} \int_{X^{k+1}} \sum_{l=1}^{l=m}
p_{\infty}(y_0,x_0)\bar{\partial}_{x_0}g_l(x_0)D_{1,x_0}Q_{\tau_1}(x_0,x_1)......$$
$$....f_l(x_k)p_{\infty}(x_k,y_0') |dx_0|....|dx_k|d\tau_1...d\tau_k $$
$$ \pm \int_{C_{k+1}} \int_{X^{k+1}} \sum_{l=1}^{l=m}
p_{\infty}(y_0,x_0)g_l(x_0)\bar{\partial}_{x_0}D_{1,x_0}Q_{\tau_1}(x_0,x_1)......$$
$$....f_l(x_k)p_{\infty}(x_k,y_0') |dx_0|....|dx_k|d\tau_1...d\tau_k $$
$$ = \pm \int_{C_{k+1}} \int_{X^{k+1}}
[\bar{\partial}^a_{x_0}p_{\infty}(y_0,x_0)][\sum_{l=1}^{l=m}f_l(x_k)g_l(x_0)][D_{1,x_0}Q_{\tau_1}(x_0,x_1)]......$$
$$....[p_{\infty}(x_k,y_0')] |dx_0|....|dx_k|d\tau_1...d\tau_k $$
$$ \pm \int_{C_{k+1}} \int_{X^{k+1}}
[p_{\infty}(y_0,x_0)][\sum_{l=1}^{l=m}f_l(x_k)g_l(x_0)][\bar{\partial}_{x_0}D_{1,x_0}Q_{\tau_1}(x_0,x_1)]......$$
\begin{equation} \label{case12}....[p_{\infty}(x_k,y_0')] |dx_0|....|dx_k|d\tau_1...d\tau_k
\end{equation}

For each of the summands dealt with by equations $\eqref{case10}$,
$\eqref{case11}$ and $\eqref{case12}$ respectively, one can argue
exactly as we did after $\eqref{case1}$. \\

 This finally proves Proposition 3.\\

\end{proof}

\begin{prop}
Let $\beta = \beta_1 +....+\beta_k$ be a $0$-cycle in
$\text{C}^{\bullet}(\diffdt(\cale))$ with $\beta_k \in
\diffdt(\cale)^{\otimes k}$. Then,

$$\fls(\beta)= \fls(\beta_1) \text{ . }$$
\end{prop}

\begin{proof}
Pick any $t >0$. If $f$ is as in Construction 1, then the cycle
$\beta$ is homotopic to the cycle $N_f(\beta)$. Since $\fls$ is a
Hochschild $0$-cocycle of $\diffdt(\cale))$,
$$\fls(\beta)=\fls(N_f(\beta))=
\fls(N_f(\beta_1))+...+\fls(N_f(\beta_k))\text{ . }$$

By proposition 4, $$\sum_{j=2}^{j=k} |\fls(N_f(\beta_j))| \leq
C\epsilon(t)$$ for some constant $C$ dependent on $\beta$ only.
Also,
$$N_f(\beta_1)= \beta_1 \pm \sum_{l=1}^{l=m}(\bar{\partial}{f_l})
\otimes g_l\beta_1 \pm \sum_{l=1}^{l=m} f_l \otimes (\bar{\partial}
g_l)\beta_1 \text{ . }$$

The proof of Proposition 4 can also be used to verify that
$$|\fls(\sum_{l=1}^{l=m}(\bar{\partial}{f_l}) \otimes g_l\beta_1 \pm
\sum_{l=1}^{l=m} f_l \otimes (\bar{\partial} g_l)\beta_1)| \leq
C'\epsilon(t) $$ for some constant $C'$ dependent only on $\beta_1$.
It follows that
$$|\fls(\beta)-\fls(\beta_1)|=|\fls(N_f(\beta)-\fls(\beta_1)| \leq
(C+C')\epsilon(t) \text{ . }$$

Now, $t$ may be chosen to be arbitrarily small. In this case,
$\epsilon(t)$ becomes arbitrarily small as well. As $C+C'$ depends
only on $\beta$, the desired proposition follows.

\end{proof}

Proposition 3 implies the following (surprising) statement.\\

\begin{prop}
The linear functional $$D \mapsto \text{str}(\Pi_0 D \Pi_0)$$
vanishes on $0$-chains of $\diffdt(\cale)$ that are commutators of
elements in $\diffdt(\cale)$.
\end{prop}

\begin{proof}
Let $D_1,D_2 \in \text{Diff}^{0}(\cale)$. This proposition is immediate from Proposition 4 and the following facts.\\

1. $\fls(d(D_1 \otimes D_2)) = 0$ since $\fls$ is
a Hochschild $0$-cocycle.\\

2. The component of $d (D_1 \otimes D_2)$ in
$\diffdt(\cale)$ is precisely $-D_1D_2+D_2D_1$.\\

\end{proof}

Of course, when $D$ is the commutator of two purely holomorphic
differential operators on $\cale$, the above proposition is
standard. The above proposition in general, is however very
counterintuitive and the author does not see any
other way of proving it.\\

\subsection{Extending the supertrace II.}

Let $\difdt(\cale)(U)$ denote $\Gamma(U,\difdt(\cale))$. Note that
the differential on $\text{C}^{\bullet}(\difdt(\cale)(U))$ extends
to a differential on the graded vector space $\oplus_k
\difdt(\cale^{\boxtimes k})(U^k)[k-1]$. We denote the resulting
complex by $\widehat{\text{C}^{\bullet}}(\difdt(\cale)(U))$. \\

Let $\widetilde{\text{hoch}(\dif(\cale))}$ denote the sheaf
associated to the presheaf
$$U \mapsto \widehat{\text{C}^{\bullet}}(\difdt(\cale)(U))$$ of
complexes of $\compl$-vector spaces. Note that
$\widetilde{\text{hoch}(\dif(\cale))}$ is a complex whose terms are
modules over the sheaf of smooth functions on $X$. It follows that
$$\text{H}^{\bullet}(\Gamma(X,\widetilde{\text{hoch}(\dif(\cale))})) \simeq
{\mathbb H}^{\bullet}(X,\widetilde{\text{hoch}(\dif(\cale))}) \text{
. }$$

The following proposition follows from Proposition 4.\\

\begin{prop} $\fls$ extends to a linear functional on
$\Gamma(X,\widetilde{\text{hoch}(\dif(\cale))})$ that vanishes on
the image of the differential $d$ of
$\Gamma(X,\widetilde{\text{hoch}(\dif(\cale))})$.
\end{prop}

\begin{proof}
\text{}\\

{\it Part 1: Constructing $\hat{\text{str}}$:}\\

Let $\widetilde{\text{hoch}^k(\dif(\cale))}$ denote the sheaf
associated to the presheaf $$U \mapsto \difdt(\cale^{\boxtimes
k})(U^k)[k-1]$$ of graded $\compl$-vector spaces. Then, any
$0$-cycle in $\Gamma(X,\widetilde{\text{hoch}(\dif(\cale))})$ is
given by $\beta_1+...+\beta_k$ where $\beta_i \in
\Gamma(X,\widetilde{\text{hoch}^i(\dif(\cale))})$. Note that
$\widetilde{\text{hoch}^1(\dif(\cale))} = \difdt(\cale)$. It follows
that $\beta_1$ is an element of $\text{Diff}^0(\cale)$. \\

Let $\hat{\text{str}}$ be the linear functional on
$\Gamma(X,\widetilde{\text{hoch}(\dif(\cale))})$ that vanishes on
$p$-chains whenever $p \neq 0$. For a zero chain $\beta =
\beta_1+...+\beta_k$ with $\beta_i \in
\Gamma(X,\widetilde{\text{hoch}^i(\dif(\cale))})$, set
$$\hat{\text{str}}(\beta)=\fls(\beta_1) \text{ . }$$

This is a well defined linear functional on the space of $0$-chains
of $\Gamma(X,\widetilde{\text{hoch}(\dif(\cale))})$.\\

{\it Part 2:}\\

We now need to show that $\hat{\text{str}}$ vanishes on the image of
$d$.\\

\textbf{Claim A:} The natural map from $\diffdt(\cale^{\boxtimes
2})[1]$ to $\Gamma(X,\widetilde{\text{hoch}^2(\dif(\cale))})$ is a
surjection of graded $\compl$-vector spaces.\\

We postpone the proof of Claim A for the next part of the proof. Let
Note that the differential on $\text{C}^{\bullet}(\diffdt(\cale))$
extends to a differential $d$ on $\oplus_k \diffdt(\cale^{\boxtimes
k})[k-1]$. This differential is a sum of the differentials
$d_{\text{H}}:\diffdt(\cale^{\boxtimes k}) \rar
\diffdt(\cale^{\boxtimes k-1})$ and $\bar{\partial}$. By Proposition
4 and Claim 1, to show that $\hat{\text{str}}$ vanishes on the image
of $d$, it suffices to show that the map $\text{Diff}^{0}(\cale)
\rar \compl$
$$D \mapsto \text{Str}(\Pi_0 D \Pi_0)$$ vanishes on the image of
$d_{\text{H}}:\text{Diff}^0(\cale^{\boxtimes 2}) \rar
\text{Diff}^0(\cale)$. \\

Let $\text{Diff}^{\leq k,\bullet}(\cale)$ denote
$\Gamma(X,\text{Dolb}^{\bullet}(X,\strc) \otimes \dif^{\leq
k}(\cale))$ where $\dif^{\leq k}(\cale)$ denotes the sheaf of
differential operators on $\cale$ of order $\leq k$. Equip
$\text{Diff}^{\leq k,\bullet}(\cale)$ with the topology generated by
the family of seminorms $\{||.||_{K,s} \text{  }| K \subset X \text{
compact }, s \in \Gamma(K,\cale \otimes \Omega^{\bullet}) \}$ given
by $$||D||_{K,s} = \text{Sup}\{||D(s)(x)|| \text{ }| x \in K\}\text{
. }$$ The topology on $\diffdt(\cale)$ is the direct limit of the
topologies on the $\text{Diff}^{\leq k, \bullet}(\cale)$. The
topology on $\diffdt(\cale^{\boxtimes 2})$ is defined
analogously.\\

Note that $\diffdt(\cale)^{\otimes 2}$ is dense in
$\diffdt(\cale^{\boxtimes 2})$. Further,
$d_{\text{H}}:\diffdt(\cale^{\boxtimes 2}) \rar \diffdt(\cale)$ is
continuous. The restriction of $d_{\text{H}}$ to
$\text{\diff}^0(\cale)^{\otimes 2}$ is just the map
$$D_1 \otimes D_2 \mapsto D_2D_1 -D_1D_2 \text{ . }$$
The desired proposition now follows from Proposition 5.\\

{\it Part 3:Proof of Claim A.}\\

Given any element $\alpha$ of
$\Gamma(X,\widetilde{\text{hoch}^2(\dif(\cale))})$, pick a finite
cover $X = \cup_i U_i$ of open sets such that $\alpha|_{U_i}$ is an
element of $\diffdt(\cale^{\boxtimes 2})(U_i \times U_i)[1]$. Pick a
partition of unity $\{f_i\}$ by compactly supported smooth functions
subordinate to he cover $\cup_i U_i$. Let $g_i$ be a compactly
supported smooth function supported on a subset of $U_i$ that is
identically $1$ on the support of $f_i$. Then,
$$\sum_i (f_i \boxtimes g_i)\alpha|_{U_i}$$ is an element of
$\diffdt(\cale^{\boxtimes 2})[1]$ whose image in
$\Gamma(X,\widetilde{\text{hoch}^2(\dif(\cale))})$ is $\alpha$.\\

\end{proof}

Next, we note that we have a natural (degree preserving) map \\
$\gamma: \widehat{\text{hoch}}(\dif(\cale)) \rar
\widetilde{\text{hoch}(\dif(\cale))}$ of complexes of sheaves of
$\compl$-vector spaces on $X$ . This induces a map $\gamma_*:
\hcoh{\bullet}{\widehat{\text{hoch}}(\dif(\cale))} \rar
\hcoh{\bullet}{\widetilde{\text{hoch}(\dif(\cale))}}$. Recall that
$\holl{0}(\Gamma(X,\widetilde{\text{hoch}(\dif(\cale))} = \hcoh{0}{\widetilde{\text{hoch}(\dif(\cale))}}$. The following corollary of Proposition 6 is immediate. \\

\begin{cor} (Corollary to Proposition 6) The Hochschild cocycle $\text{str} : \hoch{0}(\diffdt(\cale)) \rar \compl$ extends to a
$\compl$-linear functional \\
$\hat{\text{str}}:\hcoh{0}{\widetilde{\text{hoch}(\dif(\cale))}}
\rar \compl$.
\end{cor}

It is then immediate that $\hat{\text{str}} \circ \gamma_*$ gives us
a linear functional on
\\$\hcoh{0}{\widehat{\text{hoch}}(\diffdt(\cale))}$. We state this
as a corollary for emphasis.\\

\begin{cor} The linear functional $\hat{\text{str}}$ in the previous corollary gives us a linear functional $\hat{\text{tr}}:
\hcoh{0}{\widehat{\text{hoch}}(\dif(\cale))}\rar \compl$.
\end{cor}

We denote $\hcoh{-i}{\widehat{\text{hoch}}(\dif(\cale))}$ by $\choch{-i}{(\dif(\cale))}$.

Now consider a global holomorphic differential operator $D$ on
$\cale$. $D$ can be thought of as an element of $\Gamma(X,
\dif^{\bullet}(\cale))$ as well. Thus, $D$ gives us a $0$-cocycle of
the cochain complex $\Gamma(X, \prdf(\cale)^{\htens \bullet})$ . Let
$\bar{D}$ denote this cocycle. \\
One the other hand, $D$ can be thought of as a global section of the degree $0$ term of the complex $\widehat{\text{hoch}}(\dif(\cale))$. It
follows that $D$ yields a $0$-cycle in the bi-complex computing $\hcoh{\bullet}{\widehat{\text{hoch}}(\dif(\cale))}$. Call this $0$-cycle
$\tilde{D}$. It is easy to check that $\gamma(\tilde{D})= \bar{D}$. We have thus, proven the following proposition. \\

\begin{prop}

The following diagram commutes. \\

 $$\begin{CD}
\hoch{0}(\diff(\cale)) @>>> \choch{0}(\dif(\cale))\\
@VV{\text{str}}V    @V{\hat{\text{tr}}}VV \\
\compl @>>{\text{id}}> \compl \\
\end{CD}$$

\end{prop}

 \subsection{The first bookkeeping lemma}

 Consider $\cale$ as a sub bundle of $\cale \oplus {\calg F}$ where $\calg F$ is another vector bundle on $X$. We then have a map
 $\iota:\dif(\cale) \rar \dif(\cale \oplus \calg F)$ of sheaves of
 $\compl$-vector spaces whose restriction to $U$
 $\iota:\Gamma(U,\dif(\cale)) \rar \Gamma(U,\dif(\cale \oplus {\calg F}))$ is an injection that preserves addition and multiplication for each open $U \subset X$.
  This also induces an map $\bar{\iota}:\difdt(\cale) \rar \difdt(\cale \oplus \calg F)$
  of sheaves of $\compl$-vector spaces whose restriction to $U$
 $\bar{\iota}: \dif^{\bullet}(\cale)(U) \rar \dif^{\bullet}(\cale \oplus {\calg F})(U)$ is an injection for each open $U \subset X$ .\\

$\bar{\iota}$ preserves addition and multiplication though it does
not preserve the identity. $\bar{\iota}$ therefore induces a map
from $\widetilde{\text{hoch}(\dif(\cale))}$ to
$\widetilde{\text{hoch}(\dif(\cale \oplus \calg F))}$, which we will
denote by $\widetilde{\iota}$ . Let $\widetilde{\iota}_*$ denote
 the map from $\holl{*}(\Gamma(X,\widetilde{\text{hoch}(\dif(\cale))}))$ to $\holl{*}(\Gamma(X,\widetilde{\text{hoch}(\dif(\cale \oplus \calg F))}))$ induced by
 $\widetilde{\iota}$.
 The following "bookkeeping lemma"
 holds. \\

\begin{lem}
The following diagram commutes : \\

$$\begin{CD}
 \holl{0}(\Gamma(X, \widetilde{\text{hoch}(\dif(\cale))})) @> \widetilde{\iota}_* >> \holl{0}(\Gamma(X,\widetilde{\text{hoch}(\dif(\cale \oplus \calg F))} )) \\
  @VV{\hat{\text{str}}}V                                      @V{\hat{\text{str}}}VV \\
  \compl @> \id >> \compl \\
  \end{CD} $$

  \end{lem}

\begin{proof}

By Proposition 5 and the proof of Proposition 6, it is enough to
show that for any $D \in \text{Diff}^{0}(\cale)$,
$$\text{str}(\Pi_{\dol{0}{\cale}}D{\calg I}_{\dol{0}{\cale}}) = \text{str}(\Pi_{\dol{0}{\cale
\oplus \calg F}} \widetilde{\iota}(D) {\calg I}_{\dol{0}{\cale
\oplus \calg F}}) \text{ . }$$

The above equality holds since the following diagrams commute:\\

$$\begin{CD}
\dol{}{\cale} @>>> \dol{}{\cale \oplus \calg F} \\
@VV{D}V     @V{\widetilde{\iota}(D)}VV\\
\dol{}{\cale} @>>> \dol{}{\cale \oplus \calg F} \\
\end{CD}$$

$$\begin{CD}
\dol{}{\cale} @>>> \dol{}{\cale \oplus \calg F} \\
@VV{\Pi_{\dol{0}{\cale}}}V     @V{\Pi_{\dol{0}{\cale \oplus \calg F}}}VV\\
\dol{0}{\cale} @>>> \dol{0}{\cale \oplus \calg F} \\
\end{CD}$$

$$\begin{CD}
\dol{0}{\cale} @>>> \dol{0}{\cale \oplus \calg F} \\
@VV{{\calg I}_{\dol{0}{\cale}}}V     @V{{\calg I}_{\dol{0}{\cale \oplus \calg F}}}VV\\
\dol{}{\cale} @>>> \dol{}{\cale \oplus \calg F} \\
\end{CD}$$

\end{proof}

The map $\iota$ induces a map
$\hat{\iota}:\widehat{\text{hoch}}(\dif(\cale)) \rar
\widehat{\text{hoch}}(\dif(\cale \oplus \calg F)) $. Let
$\hat{\iota}_*$ denote the map induced by $\hat{\iota}$ from
$\choch{*}(\dif(\cale))$ to $\choch{*}(\dif(\cale \oplus \calg F))$.
The following corollary now follows .\\

\begin{cor}
The following diagram commutes .\\
$$ \begin{CD}
\choch{0}{(\dif(\cale))} @> \hat{\iota}_* >> \choch{0}{(\dif(\cale \oplus \calg F))} \\
@V{\hat{\text{tr}}}VV       @VV{\hat{\text{tr}}}V \\
\compl @> \id >> \compl \\
\end{CD} $$
\end{cor}

\begin{proof} For this we only need to recall that by the definition of $\iota$  at the beginning of
this section
the following diagram commutes. \\

$$\begin{CD}
 \Gamma(U^k,\dif(\cale^{\boxtimes k})) @>\hat{\iota}>> \Gamma(U^k,\dif((\cale \oplus
 \calg F)^{\boxtimes k})) \\
 @VV{\gamma}V @V{\gamma}VV \\
 \Gamma(U^k,\difdt(\cale^{\boxtimes k})) @> \widetilde{\iota} >>  \Gamma(U^k,\difdt((\cale \oplus
 \calg F)^{\boxtimes k}))
 \\
\end{CD}$$

 Taking hypercohomology, we see that $\gamma_* \circ \hat{\iota}_* =
\widetilde{\iota}_* \circ \gamma_*$. The corollary now follows
immediately from
Lemma 1. \\

\end{proof}

\section{The completed Hochschild homology of $\dif(\cale)$ and the cohomology of $X$}

\subsection{Preliminaries}

 Recall the definition of $\widehat{\text{hoch}}(\dif(\cale))$ from
 Section 3.2. Let $\widehat{\text{hoch}}(\diff(X))$ denote
 $\widehat{\text{hoch}}(\dif(\strc))$. Let $\underline{\compl}$
 denote the constant sheaf on $X$ such that
 $\Gamma(U,\underline{\compl})\simeq \compl$ for every open $U
 \subset X$ and whose restriction maps are all identity. \\

 \begin{lem}[Bryl]

 $\widehat{\text{hoch}}(\diff(X))$ is quasiisomorphic to $\underline{\compl}[2n]$ where $n$ is the dimension of $X$.

 \end{lem}

 \begin{proof} (Sketch of the Proof recalled from [Bryl]) It is enough to check the above fact for any  open ball $U \subset X$, at the level of pre sheaves.
 In other words, it is enough to show that the complex $\diff(U^{\bullet})$ with Hochschild differential is quasiisomorphic to $\compl [2n]$
  as complexes of $\compl$-vector spaces. \\

 We filter $\diff(U)$ by degree. More precisely, set $F^{-n}(\diff(U)) : = \diff^{\leq n}(U)$.The associated graded
$\text{gr}(\diff(U))$ is the space of functions
$\strcc{\text{T}^*U}$ on the cotangent bundle $\text{T}^* U$
 that are algebraic(polynomial) along the fibres.  \\

 This filtration of $\diff(U)$ is exhaustive. It yields a filtration of the complex $\diff(U^{\bullet})$ equipped with Hochschild differential.
 This filtration yields a spectral sequence converging to the cohomology of the complex $\diff(U^{\bullet})$ . $E_1^m:= \oplus_{p+q=m} E_1^{p,q}$
 term is $m$ th cohomology of the completed Hochschild complex of $\strcc{\text{T}^*U}$. This is isomorphic to the space of $-m$-holomorphic forms on
 $\text{T}^*U$ that are algebraic along the fibres (see [Bryl]). Given local coordinates $z_1,..z_n$ on $U$ and $y_1,...,y_n$ on the fibre of $\text{T}^*U$,
 defining the weight of $dz_i$ to be $0$ and that of $dy_i$ to be $1$ for all $i$ enables us to define the notion of the weight of a holomorphic form on
  $\text{T}^*U$. In the spectral sequence of this proof-sketch, $E_1^{p,m-p}$ is simply the space $-m$-forms on $\text{T}^*U$ of weight $-p$.  Therefore,
  $ E_1^{-n,-n}$ is the only nonzero summand of $E_1^{-2n}$ . By Theorem 3.1.1 and corollary 2.2.2 of [Bryl], the $E_2^m = E_{\infty}^m $ term of this
 spectral sequence is the $2n+m$ th De-Rham cohomology of $\text{T}^*U$. This is $0$ if $m \neq -2n$ and $\compl$ otherwise. This proves the
 desired lemma. Moreover, $E_2^{-n,-n}$ is the only nontrivial summand of $E_2^{-2n}$. Therefore, the $-2n$th cohomology of this complex
 can be identified with $E_2^{-n,-n}$.  \\

 \end{proof}

\begin{lem}

$\widehat{\text{hoch}}(\dif(\cale))$ is quasiisomorphic to
$\underline{\compl}[2n]$ where $n$ is the dimension of $X$.

\end{lem}

\begin{proof} This is again something that needs to be verified locally. We imitate the proof of the Morita invariance of Hochschild homology in
[Loday] Section 1.2 here. \\

{\it Part 1: Recalling the proof of Morita invariance of the
Hochschild homology of a $\compl$-algebra : } \\

 Recall from [Loday]
that if $A$ is any $\compl$-algebra, and if $\text{M}_r(A)$ denotes
the algebra of $r \times r$ matrices with entries in $A$, then we
have a map $\text{tr}$ from the Hochschild complex of
$\text{M}_r(A)$ to that of $A$. The Hochschild chain $M_1 \otimes
... \otimes M_k$ is mapped to $\text{tr} (M_1 \odot .... \odot M_k)$
where $\odot: \text{M}_r(A) \otimes \text{M}_r(B) \rar \text{M}_r(A
\otimes B)$ is an exterior multiplication. There is also a map of
complexes $\text{inc}$ in the opposite direction which is induced by
the inclusion of $A$ in $\text{M}_r(A)$ taking an element $a$ of $A$
to the matrix with $a.E_{11}$. $\text{tr} \circ \text{inc} = \id$
and there is a pre simplicial homotopy $h$ from $\text{inc} \circ
\text{tr}$ to $\id$. This homotopy is given by $h = \sum_i {(-1)}^i
h_i $ where $h_i : \text{M}_r(A)^{\otimes k+1} \rar
\text{M}_r(A)^{\otimes k+2}$ is defined by the formula
$$ h_i(\alpha^0 \otimes ... \otimes \alpha^k) = \sum E_{j1}(\alpha^0_{jr}) \otimes E_{11}(\alpha^1_{rm}) \otimes .. $$ $$.. \otimes
E_{11}(\alpha^i_{pq}) \otimes E_{1q}(1) \otimes \alpha^{i+1} \otimes  ... \otimes \alpha^k $$

The sum here is over all possible tuples of indices $(j,r,...,p,q)$.[Loday] (Section 1.2). $E_{ij}(x) = x.E_{ij}$, where $E_{ij}$ is the elementary matrix
whose only nonzero entry is a $1$ at the $ij$ position. \\

Let $U$ be an open ball contained in $X$. \\

{\it Part 2: Morita invariance of the completed Hochschild homology
of $\diff(U)$ :} \\

Recall that the completed Hochschild complex
$\widehat{\hcc{\bullet}{\diff(U)}}$ is obtained by equipping the
graded vector space $\oplus_{k \geq 1} \diff(U^k)[k-1]$ with the
Hochschild differential. In this complex, $\diff(U^k)$ should be
viewed as the $k$-th completed tensor power of $\diff(U)$. With this
in mind, the $k$-th completed tensor power of $\text{M}_r(\diff(U))$
will be $\text{M}_r(\compl)^{\otimes k} \otimes \diff(U^k)$. One may
verify that the Hochschild differential on
$\hcc{\bullet}{\text{M}_r(\diff(U))}$ extends to a differential on
the graded vector space  $\oplus_{k \geq 1}
\text{M}_r(\compl)^{\otimes k} \otimes \diff(U^k)[k-1]$. The
resulting complex is the completed Hochschild complex
$\widehat{\hcc{\bullet}{\text{M}_r(\diff(U))}}$. \\

One can verify without much difficulty that $tr$ and $inc$ extend to
maps of complexes $tr: \widehat{\hcc{\bullet}{
\text{M}_r(\diff(U))}} \rar \widehat{\hcc{\bullet}{\diff(U)}}$ and
$inc:\widehat{\hcc{\bullet}{\diff(U)}} \rar \widehat{\hcc{\bullet}{
\text{M}_r(\diff(U))}}$ respectively such that $tr \circ inc = id$.
It is also useful for us to note that explicitly, if $m_1,...,m_k
\in \text{M}_r(\compl)$,and if $\alpha \in \diff(U^k)$, then
$${tr}(m_1 \otimes ... \otimes m_k \otimes \alpha) = tr(m_1 \circ ... \circ m_k)\alpha \text{ . }$$
Moreover $h$ extends to a map $h:
\widehat{\hcc{k}{\text{M}_r(\diff(U))}} \rar
\widehat{\hcc{k+1}{\text{M}_r(\diff(U))}} $ for all $k$ with $dh +
hd = id - inc \circ tr$. Thus,\\ $tr: \widehat{\hcc{\bullet}{
\text{M}_r(\diff(U))}} \rar \widehat{\hcc{\bullet}{\diff(U)}}$  is a
quasiisomorphism. It follows from Lemma 2 that
$\widehat{\hcc{\bullet}{ \text{M}_r(\diff(U))}}$ is quasiisomorphic
to $\compl[2n]$. \\

{\it Part 3: Proof of the Lemma }\\

Let $r$ be the rank of $\cale$. If $U$ is an open ball of $X$, on
which $\cale$ is trivial, then $\diff(\cale |_U)$ is isomorphic to
$\text{M}_r(\diff(U))$ as topological algebras. The actual
isomorphism depends on the choice of (holomorphic) trivialization of
$\cale$. For a holomorphic trivialization $\phi$ of $\cale$ over
$U$, let $\phi_{\circ}$ denote the isomorphism between $\diff(\cale
|_U)$ and $\text{M}_r(\diff(U))$. Then, $\phi_{\circ}$ induces an
isomorphism of complexes $\hcc{\bullet}{\dif(\cale)(U)} \rar
\hcc{\bullet}{\text{M}_r(\diff(U))}$. This extends to an isomorphism
of complexes $\phi_*: \widehat{\hcc{\bullet}{\dif(\cale)(U)}} \rar
\widehat{\hcc{\bullet}{\text{M}_r(\diff(U))}}$. We however, know
that $\widehat{\hcc{\bullet}{\text{M}_r(\diff(U))}}$
is quasiisomorphic to $\compl[2n]$. \\

The desired lemma will follow provided that we check that the
quasiisomorphism $tr \circ \phi_*$ is independent of the choice of
trivialization $\phi$. We will show that if $\psi$ is another
holomorphic trivialization of $\cale$ over $U$, then the maps $tr
\circ \psi_* = tr \circ \phi_*$ of complexes from
$\widehat{\hcc{\bullet}{\diff(\cale |_U)}}$ to
$\widehat{\hcc{\bullet}{\diff(U)}}$ induce the same map on cohomology. \\

Note that $ \phi_* \circ \psi_*^{-1}:
\widehat{\hcc{\bullet}{\text{M}_r(\diff(U))}} \rar
\widehat{\hcc{\bullet}{\text{M}_r(\diff(U))}}$ is induced by the map
$\text{M}_r(\diff(U)) \rar \text{M}_r(\diff(U))$ taking an element
of $\text{M}_r(\diff(U))$ to its conjugate by an $r \times r$ matrix
$N$ of holomorphic functions on $U$. Denote this map of completed
Hochschild complexes by $c_N$. We only need to show that the maps
$tr \circ
c_N$ and $tr$ of complexes of $\compl$-vector spaces induce the same map on cohomology. \\

Recall that $\diff(U)$ is a filtered algebra with $F_{-m}(\diff(U))
= \diff^{\leq m}(U)$, the space of (holomorphic) differential
operators on $U$ of order at most $m$. This yields us a filtration
on $\text{M}_r(\diff(U))$ with $F_{-m} (\text{M}_r(\diff(U))) =
\text{M}_r(\diff^{\leq m}(U))$. We have
$\text{gr}(\text{M}_r(\diff(U))) =
\text{M}_r(\text{gr}(\diff(U)))$. \\

The above filtration gives us a spectral sequence with \\ $E_1^{n} =
\choch{-n}(\text{gr}(\text{M}_r(\diff(U))))$ and $E_{\infty}^n =
\choch{-n}(\text{M}_r(\diff(U)))$. The endomorphism induced by the
endomorphism $\alpha \leadsto N \alpha N^{-1}$ where $N$ is a matrix
of holomorphic functions preserves the filtration on
$\text{M}_r(\diff(U))$ and thus induces an endomorphism of a
spectral sequence on the spectral
sequence described above. \\

Note that $\text{gr}(M_r(\diff(U))) =
\text{M}_r(\text{gr}(\diff(U)))$ and that the endomorphism induced
by conjugation by $N$ is still conjugation by $N$. Denote this
endomorphism by $c_N$ as well. We now claim that the following
diagram commutes upto cohomology.

$$ \begin{CD} \widehat{\hcc{\bullet}{\text{M}_r(\text{gr}(\diff(U)))}} @>c_N >> \widehat{\hcc{\bullet}{\text{M}_r(\text{gr}(\diff(U)))}}\\
@VV{\text{tr}}V      @V{\text{tr}}VV \\
\widehat{\hcc{\bullet}{\text{gr}(\diff(U))}} @> \id >> \widehat{\hcc{\bullet}{\text{gr}(\diff(U))}}\\
\end{CD} $$

To see this, note that if $\alpha \in
\widehat{\hcc{k}{\text{gr}(\diff(U))}}$ is a cocycle, $\alpha =
tr(\frac{1}{r} id \otimes .... \otimes id \otimes \alpha)$. In other
words, $\alpha$ is the trace of a "scalar" matrix whose diagonal
elements are $\alpha$ upto a scalar factor. Now, $$c_N(id \otimes
.... \otimes id \otimes \alpha)=(N \otimes .... \otimes N) (id
\otimes .... \otimes id \otimes \alpha) (N^{-1} \otimes .... \otimes
N^{-1}) $$ $$= (id \otimes .... \otimes id \otimes \alpha) \text{ .
}$$ It follows that $(tr \circ c_N)_{*}(\overline{id \otimes ....
\otimes id \otimes \alpha}) = tr_*(\overline{id \otimes .... \otimes
id \otimes \alpha})$ where \\ $\overline{id \otimes .... \otimes id
\otimes \alpha}$ denotes the class in cohomology of the cocycle
\\ $\overline{id
\otimes .... \otimes id \otimes \alpha}$.\\

Therefore, the map of spectral sequences induced by $tr \circ c_N$
coincides with that induced by $tr$ at the $E_1$ level, and hence at
the $E_2=E_{\infty}$ level. Finally, by part 2 of this proof and by
the proof-sketch for Lemma 2, $E_2^{-n,-n}$ is the only nonzero
$E_2$ term in this spectral sequence , and $E_2^{-n,-n} \simeq
\compl$. It follows that the endomorphisms on cohomology induced by
$tr$ and $tr \circ c_N$ are indeed the endomorphisms they induce on
$E_2^{-n,-n}$. We have just shown that the maps induced $tr$ and $tr
\circ c_N$ on $E_r$ terms coincide for any $r\geq 1$.
This proves the desired lemma.\\

\end{proof}

We can now state the following immediate corollary to Lemma 3.\\

\begin{cor}

$$ \choch{-i}{(\dif(\cale))} \simeq \hh{2n-i}{\compl} \text{ . } $$

\end{cor}

We denote the isomorphism described in this corollary by
$\beta_{\cale}$.

\subsection{The second "bookkeeping" lemma}

We once more look at the situation where $\cale$ is a direct summand
of $\cale \oplus \calg F$. Notation is as in Section 3.4 of this
paper. We have a map $\iota: \dif(\cale) \rar \dif(\cale \oplus
\calg F)$. This induces a map , denoted by $\hat{\iota}$ from
$\widehat{\text{hoch}}(\dif(\cale))$ to $\widehat{\text{hoch}}(\dif(\cale \oplus \calg F))$. The following lemma holds. \\

\begin{lem}

The following diagram commutes. \\

$$ \begin{CD}
\choch{-i}{(\dif(\cale))} @>\hat{\iota}_* >> \choch{-i}{(\dif(\cale \oplus \calg F))} \\
@VV{\beta_{\cale}}V      @V{\beta_{\cale \oplus \calg F}}VV \\
\hh{2n-i}{\compl}  @> \id >>  \hh{2n-i}{\compl} \\
\end{CD} $$

\end{lem}

\begin{proof}

{\it Step 1:}\\

Let $\text{D}(\text{Sh}_{\compl}(X))$ denote the
derived category of sheaves of $\compl$-vector spaces on $X$. \\

 Recall that in the proof of Lemma 3 we
showed that the complex $\widehat{\text{hoch}}(\dif(\cale))$ of
sheaves of $\compl$ -vector spaces  was quasiisomorphic to
$\widehat{\text{hoch}}(\dif(\strc))$. Denote this quasiisomorphism
by $i_{\cale}$. $\widehat{\text{hoch}}(\dif(\strc))$ is
quasiisomorphic to
$\underline{\compl}[2n]$. Let $i$ denote this quasiisomorphism {\it for this proof}. \\

It suffices to verify that the following diagram commutes in
$\text{D}(\text{Sh}_{\compl}(X))$.

\begin{equation} \label{tochk}
\begin{CD}
\widehat{\text{hoch}}(\dif(\cale)) @>\hat{\iota}>>
\widehat{\text{hoch}}(\dif(\cale \oplus \calg F)) \\
@VV{i \circ i_{\cale}}V    @V{i \circ i_{\cale \oplus \calg F} }VV\\
\underline{\compl}[2n] @>\id >>\underline{\compl}[2n]\\
\end{CD} \end{equation}

Since a sheaf of $\compl$-vector spaces is injective iff it is
flasque (see [Riet], Lemma 3.3), the constant sheaf
$\underline{\compl}$ is an injective object in the category of
sheaves of $\compl$-vector spaces on $X$. It follows from this that
$$\text{Hom}_{\text{D}(\text{Sh}_{\compl}(X))}(\underline{\compl},\underline{\compl})
\simeq \compl \text{ . }$$ The diagram $\eqref{tochk}$ therefore,
commutes in $\text{D}(\text{Sh}_{\compl}(X))$ upto a scalar factor.
Checking that that scalar factor is one can be "done locally". It
therefore, suffices to verify that there exists a neighbourhood $U$
of every point in $X$ such that the following diagram commutes in
$\text{D}(\text{Sh}_{\compl}(U))$.

\begin{equation} \label{tochk1}
\begin{CD} \widehat{\text{hoch}}(\dif(\cale))|_U @>\hat{\iota}|_U>>
\widehat{\text{hoch}}(\dif(\cale \oplus \calg F)) |_U\\
@VV{i \circ i_{\cale} |_U}V    @V{i \circ i_{\cale \oplus \calg F} |_U}VV\\
\underline{\compl}[2n] @>\id >> \underline{\compl}[2n]\\
\end{CD} \end{equation}

{\it Step 2: Verifying $\eqref{tochk1}$.}\\

It suffices to verify $\eqref{tochk1}$ at the level of pre-sheaves.
We must therefore , prove that the following diagram commutes in the
category of complexes of $\compl$-vector spaces upto cohomology.

\begin{equation}
\begin{CD} \widehat{\hcc{\bullet}{\dif(\cale)(U)}}
@>\hat{\iota} >>\widehat{\hcc{\bullet}{\dif(\cale \oplus \calg F)(U)}} \\
@VV{i \circ i_{\cale} |_U}V     @V{i \circ i_{\cale \oplus \calg F} |_U}VV\\
\compl[2n]  @>\id>> \compl[2n]\\
\end{CD} \end{equation}

To verify that the above diagram commutes upto cohomology , it
suffices to verify that the diagram below commutes upto cohomology
in the category of complexes of $\compl$-vector spaces.\\

\begin{equation} \label{tochk2}
\begin{CD} \widehat{\hcc{\bullet}{\dif(\cale)(U)}}
@>\hat{\iota} >>\widehat{\hcc{\bullet}{\dif(\cale \oplus \calg F)(U)}} \\
@VV{ i_{\cale} |_U}V     @V{i_{\cale \oplus \calg F} |_U}VV\\
\widehat{\hcc{\bullet}{\dif(U)}}  @>\id>> \widehat{\hcc{\bullet}{\dif(U)}}\\
\end{CD} \end{equation}

 Most of the hard work
necessary for this step has been done already. Let $\phi_{\cale}$
and $\phi_{\calg F}$ denote holomorphic trivializations over $U$ of
$\cale$ and $\calg F$ respectively. Then, $\phi:=\phi_{\cale} \oplus
\phi_{\calg F}$ is a holomorphic
trivialization of $\cale \oplus \calg F$ over $U$. \\

This yields the following commutative diagram, all of whose
morphisms are continuous. $r$ and $s$ denote the ranks of $\cale$
and $\calg F$ respectively. \\

$$ \begin{CD}
\dif(\cale)(U) @> \iota >> \dif(\cale \oplus \calg F)(U) \\
@VV{\phi_{\cale}}V           @V{\phi}VV \\
 \text{M}_{r}(\diff(U)) @>\iota_{r,s} >> \text{M}_{r+s}(\diff(U)) \\
\end{CD} $$

This yields the following commutative diagram. \\

$$ \begin{CD}
\widehat{\hcc{\bullet}{\dif(\cale)(U)}} @>\hat{\iota}_* >>
\widehat{\hcc{\bullet}{\dif(\cale \oplus \calg F)(U)}} \\
@VV{\phi_{\cale,*}}V           @V{\phi_*}VV \\
\widehat{\hcc{\bullet}{\text{M}_r(\diff(U))}}
@>{\hat{\iota}_{r,s,*}}>>
\widehat{\hcc{\bullet}{\text{M}_{r+s}(\diff(U))}} \\
\end{CD} $$

$\iota_{r,s}$ denotes the embedding from $\text{M}_{r}(\diff(U))$ to
$\text{M}_{r+s}(\diff(U))$ which takes a matrix $\alpha \in
\text{M}_{r}(\diff(U))$ to $\bar{\alpha} \in
\text{M}_{r+s}(\diff(U))$ where ${\bar{\alpha}}_{ij} = \alpha_{ij}$
if
$(i,j) \in \{ 1,...,r\} \times \{ 1,..., r\} $ and ${\bar{\alpha}}_{ij} =0 $ otherwise. \\

Now, we proved that $i_{\cale} |_U = tr_r \circ \phi_{\cale,*}$ in
part 3 of the proof Lemma 3. Here, $tr_r$ is the map $tr:
\widehat{\hcc{\bullet}{\text{M}_{r}(\diff(U))}} \rar
\widehat{\hcc{\bullet}{\diff(U)}}$ described in the proof of Lemma
3. Similarly, $i_{(\cale \oplus \calg F) |_U} = tr_{r+s} \circ
\phi_*$.
\\

It therefore, suffices to prove that the following diagram commutes.
\\

$$ \begin{CD}
\widehat{\hcc{\bullet}{\text{M}_r(\diff(U))}}
@>{\hat{\iota}_{r,s,*}}>>
\widehat{\hcc{\bullet}{\text{M}_{r+s}(\diff(U))}} \\
@VV{tr_r}V    @V{tr_{r+s}}VV \\
\widehat{\hcc{\bullet}{\diff(U)}} @>id>>
\widehat{\hcc{\bullet}{\diff(U)}} \\
\end{CD} $$

This follows immediately from the explicit formula for $tr_r$
recalled in Part 2 of the proof of Lemma 3. This finally verifies
$\eqref{tochk2}$,
thus proving the desired lemma.\\

\end{proof}

\subsection{Proof of theorem 2}

\begin{proof}(Proof of theorem 2)

Corollary 5 states that the following diagram commutes. \\

$$\begin{CD}
\choch{0}{(\dif(\cale))} @> \hat{\iota}_* >> \choch{0}{(\dif(\cale \oplus \calg F))} \\
@VV{\hat{\text{tr}}}V               @V{\hat{\text{tr}}}VV \\
\compl @> \id >> \compl \\
\end{CD}  $$

Lemma 4 for $i= 0$ now states that the following diagram commutes. \\

$$ \begin{CD}
\choch{0}{(\dif(\cale))} @> \hat{\iota}_* >> \choch{0}{(\dif(\cale \oplus \calg F))} \\
@VV{\beta_{\cale}}V      @V{\beta_{\cale \oplus \calg F}}VV \\
\hh{2n}{\compl}  @> \id >>  \hh{2n}{\compl} \\
\end{CD} $$

It follows that $\hat{\text{tr}} \circ {\beta_{\cale}}^{-1} =
\hat{\text{tr}} \circ {\beta_{\cale \oplus \calg F}}^{-1} $. In the
notation of Theorem 2, this says that $I_{\cale} = I_{\cale \oplus
\calg F}$. A symmetric argument shows that $I_{\calg F} = I_{\cale
\oplus \calg F}$. This completes the proof of Theorem 2.

\end{proof}

\section{The completed Cyclic homology of $\diff(\cale)$}

\subsection{Recollections on Cyclic homology}

\subsubsection{Tsygan's double complex}

Let $\calg A$ be a dg-$\compl$-algebra such that ${\calg A}^n = 0 $
for almost all $n$. Let $d_{\text{bar}}$ and $d_{\text{hoch}}$
denote the differentials of $\text{bar}^{\bullet}(\calg A)$ and
$\hcc{\bullet}{\calg A}$
respectively. \\

Let $\tau: {\calg A}^{\otimes k} \rar {\calg A}^{\otimes k}$ denote
the map \\ $a_1 \otimes ....\otimes a_k \leadsto
{(-1)}^{(d_k+1)(d_1+ ... +d_{k-1} + k-1)} a_k \otimes a_1 \otimes
.... \otimes a_{k-1} $ for homogenous elements $a_1,...,a_k$ of
$\calg A$ of degrees $d_1,..,d_k$ respectively. Let $N: {\calg
A}^{\otimes k} \rar {\calg
A}^{\otimes k}$ be the map $N = 1+ \tau + .. + \tau^{k-1}$. \\

$\tau$ and $N$ induce maps from ${\calg A}^{\otimes k}$ to ${\calg
A}^{\otimes k}$ for any $k$. It follows that $\tau$ and $N$ induce
maps from $\text{bar}^{-n}(\calg A)$ to $\hcc{-n}{\calg A}$ and
$\hcc{-n}{\calg A}$ to $\text{bar}^{-n}(\calg A)$ respectively for any $n$.\\

Consider the following double complex the degree of whose non-zero
columns is non positive.

$$ \begin{CD}
..... @>>> ..... @>>> .... @>>> .... @>>> ....\\
 @A{d_{\text{hoch}}}AA   @A{-d_{\text{bar}}}AA  @A{d_{\text{hoch}}}AA
 @A{-d_{\text{bar}}}AA  @A{d_{\text{hoch}}}AA \\
.... @> N>> \text{bar}^{0}(\calg A) @>id - \tau >> \hcc{0}{\calg A}
@>N>> \text{bar}^{0}(\calg A) @>id -\tau>> \hcc{0}{\calg A} \\
@A{d_{\text{hoch}}}AA  @A{-d_{\text{bar}}}AA @A{d_{\text{hoch}}}AA
 @A{-d_{\text{bar}}}AA @A{d_{\text{hoch}}}AA \\
.... @> N>> \text{bar}^{-1}(\calg A) @>id - \tau >> \hcc{-1}{\calg
A}
@>N>> \text{bar}^{-1}(\calg A) @>id -\tau>> \hcc{-1}{\calg A} \\
@A{d_{\text{hoch}}}AA  @A{-d_{\text{bar}}}AA @A{d_{\text{hoch}}}AA
 @A{-d_{\text{bar}}}AA @A{d_{\text{hoch}}}AA \\
.... @> N>> \text{bar}^{-2}(\calg A) @>id - \tau >> \hcc{-2}{\calg
A}
@>N>> \text{bar}^{-2}(\calg A) @>id -\tau>> \hcc{-2}{\calg A} \\
@A{d_{\text{hoch}}}AA  @A{-d_{\text{bar}}}AA @A{d_{\text{hoch}}}AA
 @A{-d_{\text{bar}}}AA @A{d_{\text{hoch}}}AA \\
..... @>>> ..... @>>> .... @>>> .... @>>> .... \\
 \end{CD} $$

We denote the above double complex by
$\text{CC}^{\bullet,\bullet}(\calg A)$. $\text{CC}^{p,q}(\calg A) =
\hcc{q}{\calg A}$ if $p$ is even and non-positive,
$\text{CC}^{p,q}(\calg A) = \text{bar}^{q}(\calg A)$ if $p$ is odd
and negative, and $\text{CC}^{p,q}(\calg A) = 0$ otherwise. \\

 The vertical differential $d_v: \text{CC}^{p,q}(\calg
A) \rar \text{CC}^{p,q+1}(\calg A)$ is $d_{\text{hoch}}$ if $p$ is
even and non positive , and $-d_{\text{bar}}$ if $p$ is odd and
negative. \\

The horizontal differential $d_h: \text{CC}^{p,q}(\calg A) \rar
\text{CC}^{p+1,q}(\calg A)$ is given by $id - \tau $ if $p$ is odd
and negative and $N$ if $p$ is even and negative. \\

The double complex $\text{CC}^{\bullet,\bullet}(\calg A)$ is called
the {\it Tsygan's double complex} of $\calg A$. \\

{\it Definition } The {\it Cyclic complex} of ${\calg A}$ is the
total complex of $\text{CC}^{\bullet,\bullet}(\calg A)$. It is
denoted by $\text{Cycl}^{\bullet}(\calg A)$.

{\it Definition } The cyclic homology $\cyc{-i}(\calg A)$ is the
$-i$th cohomology of $\text{Cycl}^{\bullet}(\calg A)$. As in Section
2.1, contrary to the standard practice , we refer to a $-i$ cocycle
of $\text{Cycl}^{\bullet}(\calg A)$ as an $i$-cycle in
$\text{Cycl}^{\bullet}(\calg A)$. The following standard
 propositions are important to us.  \\

\begin{prop}

(i) If $\calg B = \enn(V^{\bullet})$ where $V^{\bullet}$ is a finite
dimensional graded $\compl$- vector space with zero differential,
then,
$$\cyc{-2i}(\calg B) \simeq \compl $$ for
all $i \geq 0$. All other cyclic homologies of $\calg B$ vanish. \\

 (ii) Further, the map $\text{tr}_{2i}$ yielding the above isomorphism, is obtained by mapping
 the class in $\cyc{-2i}(\calg B)$ of a $2i$-cycle in $\text{Cycl}^{\bullet}(\calg B)$ given by
 a tuple $(b_{2i} , b_{2i-1} ,...., b_0)$ to $\text{str}(b_0)$ where $b_k \in \text{CC}^{-2i+k,-k}(\calg
 B)$.\\

 (iii) Still further, the map $\text{tr}_{2i}$ yielding the above isomorphism,
 maps
 the class in $\cyc{-2i}(\calg B)$ of a $2i$-cycle in $\text{Cycl}^{\bullet}(\calg B)$ given by
 a tuple \\ $(b_{2i} , b_{2i-1} ,...., b_0,b_{-1},...,b_{-2j})$ ($j>0 \text{ , } b_k \in \text{CC}^{-2i+k,-k}(\calg
 B)$) to
 $\text{str}(b_0)$.

\end{prop}

\begin{proof}
There is a spectral sequence converging to $\cyc{\bullet}(\calg B)$
such that $E_1^{p,q} = {\calg H}^q(\text{CC}^{p, \bullet}(\calg B))$
. This is the spectral sequence that arises out of the filtration of
the double complex $\text{CC}^{\bullet,\bullet}(\calg B)$ by
columns. Now , since $\calg B$ has a unit
$\text{bar}^{\bullet}(\calg B)$ is acyclic by recollection 1 of
Section 2.1. $\hcc{\bullet}{\calg B}$ is quasiisomorphic to $\compl$
concentrated in degree $0$ by Proposition 1. It follows that
$E_1^{-2i,0} \simeq \compl$ for all $i \geq 0$ and $E_1^{p,q} = 0$
for all other $(p,q)$. This spectral sequence therefore collapses at
$E_1$. It follows that $\cyc{-2i}(\calg B) \simeq E_1^{-2i,0} \simeq
\compl$ for all
$i \geq 0$. This proves part (i).\\

Let $F^{\bullet}$ be the filtration on $\text{Cycl}^{\bullet}(\calg
B)$ yielding the spectral sequence in the proof of part (i) of this
proposition. Consider a tuple $(b_{2i} , b_{2i-1} ,...., b_0)$ with
$b_k \in \text{CC}^{-2i+k,-k}(\calg B)$ that yields a cyclic cycle.
Then, $(b_{2i} , b_{2i-1} ,...., b_0) \in
F^{-2i}\text{Cycl}^{\bullet}(\calg B)$. It follows that the image of
$(b_{2i} , b_{2i-1} ,...., b_0)$ in
$\frac{F^{-2i}\text{Cycl}^{\bullet}(\calg
B)}{F^{-2i+1}\text{Cycl}^{\bullet}(\calg B)} \simeq
\text{C}^{\bullet}(\calg B)$ is the Hochschild $0$-cycle $b_0$. Note
that $E_1^{-2i+j,-j} = 0$ for $j \neq 0$ and $E_1^{-2i,0} =
\text{H}^0(\frac{F^{-2i}\text{Cycl}^{\bullet}(\calg
B)}{F^{-2i+1}\text{Cycl}^{\bullet}(\calg B)}) \simeq \hoch{0}(\calg
B)$. It follows that the image of the tuple $(b_{2i} , b_{2i-1}
,...., b_0)$ in $E_1^{-2i,0}$, and therefore in $\cyc{-2i}(\calg B)$
is the image of $b_0$ in $\hoch{0}(\calg B)$. By proposition 1, this
is precisely $\text{str}(b_0)$. This proves that
$tr_{2i}(\widetilde{(b_{2i},...,b_0)})= \text{str}(b_0)$ where
$\widetilde{(b_{2i},.,b_0)}$ is the class in $\cyc{-2i}(\calg B)$ of
the cycle obtained from $(b_{2i} , b_{2i-1}
,...., b_0)$. This proves part (ii).\\

To prove part (iii), note that if $(b_{2i} , b_{2i-1} ,....,
b_0,b_{-1},...,b_{-2j})$ ($j >0$) is a cyclic cycle, then $b_{-2j}$
is a Hochschild $-2j$-cycle. It follows from proposition 1 that
$b_{-2j} = d_{\text{hoch}}c_{-2j+1} $ for some $c_{-2j+1} \in
\text{C}^{2j-1}(\calg B)$. Consider $c_{-2j+1} \in
\text{CC}^{-2i-2j,2j -1}(\calg B)$ as an element of
$\text{Cycl}^{-2i-1}(\calg B)$. Then, the cycle $(b_{2i} , b_{2i-1}
,...., b_0,b_{-1},...,b_{-2j})- d_{\text{cycl}}c_{-2j+1}$ arises out
of the cycle \\ $(b_{2i} , b_{2i-1} ,...., b_0,b_{-1},...,b_{-2j+1}+
Nc_{-2j+1},0)$ of $ \text{CC}^{\bullet,\bullet}(\calg B)$.
Therefore, the class of $(b_{2i} , b_{2i-1} ,....,
b_0,b_{-1},...,b_{-2j})$ in $\cyc{-2i}(\calg B)$ is the same as the
class of $(b_{2i} , b_{2i-1} ,...., b_0,b_{-1},...,b_{-2j+1}+
Nc_{-2j+1},0)$ in $\cyc{-2i}(\calg B)$. Since the bar complex of
$\calg B$ is acyclic , $b_{-2j+1}+
Nc_{-2j+1}=d_{\text{bar}}c_{-2j+2}$ for some element $c_{-2j+2}$ of
$\text{bar}^{2j-2}(\calg B)$. The previous step can be repeated to
show that the class of $(b_{2i} , b_{2i-1} ,....,
b_0,b_{-1},...,b_{-2j})$ in $\cyc{-2i}(\calg B)$ is the same as the
class of $(b_{2i} , b_{2i-1} ,...., b_0,b_{-1},...,b_{-2j+2}+ (\id
-\tau) c_{-2j+2},0,0)$ in $\cyc{-2i}(\calg B)$. This process can be
continued to show that the class of $(b_{2i} , b_{2i-1} ,....,
b_0,b_{-1},...,b_{-2j})$ in $\cyc{-2i}(\calg B)$ is the same as the
class of $(b_{2i} , b_{2i-1} ,...., b_0 + (\id -\tau)c_0)$ in
$\cyc{-2i}(\calg B)$ for some element $c_0$ of $\text{bar}^0(\calg
B)$. By part (ii) of this proposition, this is equal to
$\text{str}(b_0+(\id -\tau)c_0)$. But $\text{str}((\id-\tau)c_0)=0$
for any element $c_0$ of $\text{bar}^0(\calg B)$. This proves
(iii).\\

\end{proof}

If $\calg A$ is a (unital) dg-$\compl$-algebra , let  $\text{CC}^{
\{ 2\},\bullet}(\calg A)$ denote the bicomplex consisting of the
columns $\text{CC}^{-1,\bullet}(\calg A)$ and
$\text{CC}^{0,\bullet}(\calg A)$ with differentials as in
$\text{CC}^{\bullet,\bullet}(\calg A)$.
\\

 We recall that we have an exact sequence of complexes
$$
\begin{CD} 0 @>>> \text{Tot}(\text{CC}^{ \{
2\},\bullet}(\calg A))  @>I>> \text{Cycl}^{\bullet}(\calg A) @> {S}
>> \text{Cycl}^{\bullet + 2}(\calg A) @>>> 0 \end{CD} $$

We recall that $\text{CC}^{1,\bullet}(\calg A)$ is acyclic
(Recollection 1, Section 2.1). It follows that $\text{CC}^{ \{
2\},\bullet}(\calg A)$ is quasiisomorphic to $\hcc{\bullet}{\calg
A}$. This quasiisomorphism is realized by the map of complexes
taking $a_k \in \hcc{-k}{\calg A}$ to $(a_k,0) \in \hcc{-k}{\calg A}
\oplus \text{bar}^{1-k}(\calg A)$ . We denote the composite of $I$
with this quasiisomorphism by ${\calg I}$. The map $S$ is obtained
by projection to the double complex obtained from
$\text{CC}^{\bullet,\bullet}(\calg A)$ by truncating the columns
$\text{CC}^{i,\bullet}(\calg A)$ for $i = 0,-1$.
We now obtain the following proposition. \\

\begin{prop}

If $\calg B = \enn(V^{\bullet})$ where $V^{\bullet}$ is a finite
dimensional graded $\compl$- vector space with zero differential,
then, the following diagrams commute. \\

$$\begin{CD} \cyc{-2i}(\calg B) @>{S} >> \cyc{-2i +2}(\calg B) \\
           @VV{\text{tr}_{2i}}V        @V{\text{tr}_{2i-2}}VV \\
           \compl  @> \id >> \compl \\
           \end{CD} \\
           $$

$$\begin{CD} \hoch{0}(\calg B) @>{{\calg I}}>> \cyc{0}(\calg B) \\
          @VV{\text{str}}V     @V{\text{tr}_0}VV \\
\compl  @> \id >> \compl \\
           \end{CD} \\
           $$
    \end{prop}

\begin{proof}
Let $(b_{2i} , b_{2i-1} ,...., b_0)$ , $b_j \in
\text{CC}^{-2i+j,-j}(\calg B)$ be a tuple yielding a cyclic cycle.
We already demonstrated while proving part (iii) of proposition 8
that any class in $\cyc{-2i}(\calg B)$ can be represented by a cycle
coming from a tuple of this form.
 Let \\
$\widetilde{(b_{2i} , b_{2i-1} ,...., b_0)}$ denote the class of
$(b_{2i} , b_{2i-1} ,...., b_0)$ in $\cyc{-2i}(\calg B)$. Then,
$\text{tr}_{2i}(\widetilde{(b_{2i} , b_{2i-1} ,...., b_0)}) =
\text{str}(b_0)$ by Proposition 8. \\ $S((b_{2i} , b_{2i-1} ,....,
b_0)) =(b_{2i-2} , b_{2i-3} ,...., b_0)$ by the definition of $S$.
\\
$\text{tr}_{2i-2}(\widetilde{(b_{2i-2} , b_{2i-3} ,...., b_0)}) =
\text{str}(b_0)$ by Proposition 8. This proves that the first diagram commutes.\\

Let $b_0$ be a Hochschild $0$-cycle of $\calg B$. Then, ${\calg
I}(b_0) = b_0 \in \text{CC}^{0,0}(\calg B)$. Let $\bar{b_0}$ denote
the class of $b_0$ in $\hoch{0}(\calg B)$. Then, $\bar{b_0} =
\text{str}(b_0)$. However, ${\calg I}(\bar{b_0}) = \widetilde{b_0}$.
Now $\text{tr}_0(\widetilde{b_0}) =
\text{str}(b_0) $ by Proposition 8. This proves that the second diagram commutes. \\

\end{proof}

 Further, if $\calg F: \calg A \rar \calg B$ is an $A_{\infty}$ morphism with Taylor components ${\calg F}_k$, we have the following proposition. \\

 \begin{prop}

 The map ${\calg F}_{\text{Hoch}}$ mentioned in Proposition 2 extends to a map ${\calg F}_{\text{cycl}}$ of complexes from $\text{Cycl}^{\bullet}(\calg A)$
 to $\text{Cycl}^{\bullet}(\calg B)$ \\

 \end{prop}

 \begin{proof}

It suffices to check that ${\calg F}_{\text{Hoch}}$ extends to a map
of ${\calg F}_{\text{tsyg}}$ bicomplexes from
$\text{CC}^{\bullet,\bullet}(\calg A)$ to
$\text{CC}^{\bullet,\bullet}(\calg B)$. \\

Let ${\calg F}_{\text{tsyg}}: \text{CC}^{p,q}(\calg A) \rar
\text{CC}^{p,q}(\calg B)$  be ${\calg F}_{\text{Hoch}}$ if $p$ is
even and non positive and ${\calg F}_{\text{bar}}$ if $p$ is odd and
negative. \\

By Proposition 2, ${\calg F}_{\text{Hoch}}:
\text{CC}^{p,\bullet}(\calg A) \rar \text{CC}^{p,\bullet}(\calg B)$
is a map of complexes of $\compl$-vector spaces if $p$ is even and
non positive. By the definition of an $\ainf$-morphism, ${\calg
F}_{\text{bar}}: \text{CC}^{p,\bullet}(\calg A) \rar
\text{CC}^{p,\bullet}(\calg B)$ is a map of complexes if $p$ is odd
and negative. \\

The following verifications, which we leave to the reader, complete
the proof that ${\calg F}_{\text{tsyg}}$ is a map of bi-complexes,
and thus yields a map ${\calg F}_{\text{cycl}}:
\text{Cycl}^{\bullet}(\calg A) \rar \text{Cycl}^{\bullet}(\calg B)$
of complexes.

 $$ (i) \text{  } (id- \tau) \circ {\calg F}_{\text{bar}} = {\calg
F}_{\text{Hoch}} \circ (1- \tau) \text{ . }$$   $$ (ii) \text{  } N
\circ {\calg F}_{\text{Hoch}} = {\calg F}_{\text{Bar}} \circ N
\text{ . }$$ Here, $N$ and $\tau$ are as in the definitions of
$\text{CC}^{\bullet,\bullet}(\calg A)$ and
$\text{CC}^{\bullet,\bullet}(\calg B)$.

\end{proof}

We also state the following consequence of Propositions 10 and 8 as
a proposition. \\

\begin{prop} Let $\calg A$  be a dg-$\compl$ algebra. Let $\calg B$ be as in Proposition 9.  Suppose that $\calg F$ is an $A_{\infty}$ morphism from $\calg A$ to $\calg
B$. Let $\widetilde{{\calg F}_{\text{cycl}}}$ denote the map from $\cyc{\bullet}(\calg A)$ to $\cyc{\bullet}(\calg B)$ induced by
${\calg F}_{\text{cycl}}$.\\

1.   $$ \text{   } \text{tr}_{2i} \circ \widetilde{{\calg
F}_{\text{cycl}}} (\widetilde{(a_{2i}, ....,
a_0,a_{-1},...,a_{-l})}) = \text{str} ({\calg F}_{\text{Hoch}}
(a_0)) \text{ . }$$ 2. If $a_0$ is a Hochschild $0$- cycle arising
out of a degree $k-1$ element of ${\calg A}^{\otimes k}$ then,
$$ \text{tr}_{2i} \circ \widetilde{{\calg F}_{\text{cycl}}}
(\widetilde{(a_{2i},..., a_0,a_{-1},...,a_{-l})}) =
\sum_{s=0}^{s=k-1} \text{str}({\calg F}_k(\tau^s(a_0)))$$ where
$\tau$ is as in
Corollary 2.\\

\end{prop}

\begin{proof} By the proof of Proposition 10, $${\calg F}_{\text{tsyg}}(a_{2i},....,a_0,a_{-1},...,a_{-l})
= $$ $$ ({\calg F}_{\text{Hoch}}(a_{2i}),{\calg
F}_{\text{bar}}(a_{2i-1}),...,{\calg F}_{\text{bar}}(a_1),{\calg
F}_{\text{Hoch}}(a_0),....,{\calg F}_{\text{bar/hoch}}(a_{-l}))$$
where ${\calg F}_{\text{bar/hoch}}(a_{-l}) = {\calg
F}_{\text{hoch}}(a_{-l})$ if $l$ is even and  ${\calg
F}_{\text{bar/hoch}}(a_{-l}) = {\calg F}_{\text{bar}}(a_{-l})$ if
$l$ is odd. Part 1 of this proposition now follows immediately from
Proposition 8, part (iii). Part 2 of this proposition is
immediate from Part 1 and Corollary 2 . \\

\end{proof}

\subsection{The complex $\widetilde{\text{Cycl}^{\bullet}(\dif(\cale))} $}

Let $\underline{\text{C}}^{\infty}$ denote the sheaf of smooth
functions on $X$. For each open $U \subset X$, we can consider the
completed Hochschild complex of $\dif^{\bullet}(\cale)(U)$ as in
Section 3.3. Denote this complex by
$\widehat{\hcc{\bullet}{\dif^{\bullet}(\cale)}(U)}$. Consider the
sheafification of the presheaf $U \leadsto
\widehat{\hcc{\bullet}{\dif^{\bullet}(\cale)}(U)}$ .  Recall that
this sheaf of complexes of $\underline{\text{C}}^{\infty}$-modules
was denoted by $\widetilde{\text{hoch}(\dif(\cale))}$ in Section
3.3. {\it Unlike in Section 3.3}, let
$\widetilde{\text{hoch}^n(\dif(\cale))}$ denote the degree $n$
component of $\widetilde{\text{hoch}(\dif(\cale))}$.

One can also consider the completed Bar complex
$\widehat{\text{bar}^{\bullet}(\dif^{\bullet}(\cale)(U))}$. This is
defined as in Section 3.3 - the underlying graded $\compl$-vector
space of the complex
$\widehat{\text{bar}^{\bullet}(\dif^{\bullet}(\cale)(U))}$ is the
same as that of $\widehat{\hcc{\bullet}{\dif^{\bullet}(\cale)}(U)}$
but the differential on
$\widehat{\text{bar}^{\bullet}(\dif^{\bullet}(\cale)(U))}$ is the
bar differential. We will denote the sheafification of the presheaf
$U \leadsto
\widehat{\text{bar}^{\bullet}(\dif^{\bullet}(\cale)(U))}$ by
$\widetilde{\text{bar}^{\bullet}(\dif(\cale))}$.
 We have the following proposition.  \\

\begin{prop}

The differentials in Tsygan's double complex for
$\dif^{\bullet}(\cale)(U)$ extend to yield differentials for the
following double complex of $\underline{\text{C}}^{\infty}$ -modules
on
$X$. \\

$$\begin{CD}
..... @>>> ..... @>>> .... @>>> .... @>>> ....\\
 @A{d_{\text{hoch}}}AA   @A{-d_{\text{bar}}}AA  @A{d_{\text{hoch}}}AA
 @A{-d_{\text{bar}}}AA  @A{d_{\text{hoch}}}AA \\
.... @> N>> \widetilde{\text{bar}^{0}(\dif(\cale))} @>id - \tau >>
\widetilde{\text{hoch}^{0}(\dif(\cale))}
@>N>> \widetilde{\text{bar}^{0}(\dif(\cale))} @>id -\tau>> \widetilde{\text{hoch}^{0}(\dif(\cale))}  \\
@A{d_{\text{hoch}}}AA  @A{-d_{\text{bar}}}AA @A{d_{\text{hoch}}}AA
 @A{-d_{\text{bar}}}AA @A{d_{\text{hoch}}}AA \\
.... @> N>> \widetilde{\text{bar}^{-1}(\dif(\cale))} @>id - \tau >>
\widetilde{\text{hoch}^{-1}(\dif(\cale))}
@>N>> \widetilde{\text{bar}^{-1}(\dif(\cale))} @>id -\tau>> \widetilde{\text{hoch}^{-1}(\dif(\cale))}  \\
@A{d_{\text{hoch}}}AA  @A{-d_{\text{bar}}}AA @A{d_{\text{hoch}}}AA
 @A{-d_{\text{bar}}}AA @A{d_{\text{hoch}}}AA \\
.... @> N>> \widetilde{\text{bar}^{-2}(\dif(\cale))} @>id - \tau >>
\widetilde{\text{hoch}^{-2}(\dif(\cale))}
@>N>> \widetilde{\text{bar}^{-2}(\dif(\cale))} @>id -\tau>> \widetilde{\text{hoch}^{-2}(\dif(\cale))}  \\
@A{d_{\text{hoch}}}AA  @A{-d_{\text{bar}}}AA @A{d_{\text{hoch}}}AA
 @A{-d_{\text{bar}}}AA @A{d_{\text{hoch}}}AA \\
..... @>>> ..... @>>> .... @>>> .... @>>> .... \\
    \end{CD}\\ $$

 \end{prop}

\begin{proof} This only needs to be checked at the level of
presheaves. The double complex mentioned in this proposition is the
sheafification of the presheaf $U \leadsto
\widehat{\text{CC}^{\bullet,\bullet}(\dif^{\bullet}(\cale)(U))}$.
Here,
$\widehat{\text{CC}^{\bullet,\bullet}(\dif^{\bullet}(\cale)(U))}$ is
the double complex obtained by replacing the Hochschild and bar
complexes that make up the columns of Tsygan's double complex of
$\dif^{\bullet}(\cale)(U)$ with their completed versions. It is also
clear that $id-\tau$ and $N$ extend to horizontal differentials on
the double complex in this proposition.
\end{proof}

We denote the total complex of this double complex by
$\widetilde{\text{Cycl}}^{\bullet}(\dif(\cale))$. Since it is a
complex of $\underline{\text{C}}^{\infty}$ modules, its
Hypercohomology is computed by
 the complex $\Gamma(X, \widetilde{\text{Cycl}}^{\bullet}(\dif(\cale)))$. \\

The following cyclic analog of Proposition 6 follows from
Proposition 11 and the fact that $\sum_s \tau^s (id-\tau) = 0$. We
denote the $-2i$ th cohomology
of $\Gamma(X, \widetilde{\text{cycl}}^{\bullet}(\dif(\cale)))$ by $\ccyc{-2i}(\dif(\cale))$. \\

\begin{prop}
The formula described in Proposition 11  extends to yield us a
$\compl$-linear functional $\tilde{\text{tr}}_{2i}:
\ccyc{-2i}(\dif(\cale)) \rar \compl$.
\end{prop}

\begin{proof} The fact that the formula makes sense follows directly from Proposition 11 and Proposition 6.
To show that it vanishes on coboundaries, we recall that Proposition
6 also tells us that it vanishes on the image of the
$d_{\text{hoch}}$ differential. We only need to verify that it
vanishes on the image of the $id-\tau$ differential. This is a
consequence of the fact
that $\sum_s \tau^s (id-\tau) = 0$.\\
\end{proof}

Recall that we had a map $S: \text{Cycl}^{\bullet}(\diffdt(\cale))
\rar \text{Cycl}^{\bullet +2}(\diffdt(\cale))$. It is easy to verify
that this yields us a map $S:
\widetilde{\text{Cycl}}^{\bullet}(\diffdt(\cale)) \rar
\widetilde{\text{Cycl}}^{\bullet+2}(\diffdt(\cale)) $ . Similarly,
${\calg I}$ can be seen to extend to a map of complexes ${\calg I}:
\widetilde{\text{hoch}(\dif(\cale))} \rar
\widetilde{\text{Cycl}}^{\bullet}(\diffdt(\cale))$. A direct
consequence of the formula in Proposition 11 is the following proposition. \\

\begin{prop}
The following diagrams commute. \\
$$ \begin{CD}
\ccyc{-2i}(\dif(\cale)) @> {S}>> \ccyc{-2i+2}(\dif(\cale)) \\
@VV{\tilde{\text{tr}}_{2i}}V      @V{\tilde{\text{tr}}_{2i-2}}VV \\
\compl  @> \id >> \compl \\
\end{CD} $$

$$\begin{CD}
\holl{0}(\Gamma(X,\widetilde{\text{hoch}(\dif(\cale))})) @>{\calg
I}>> \ccyc{0}(\dif(\cale)) \\
@VV{\hat{\text{str}}}V @V{\tilde{\text{tr}_0}}VV \\
\compl  @> \id >> \compl \\
\end{CD} $$

\end{prop}

\subsection{The completed cyclic homology of $\dif(\cale)$ and the cohomology of $X$}

\subsubsection{The completed cyclic homology of $\dif(\cale)$}

We can define the completed Bar complex of $\dif(\cale)$. This is
the sheafification of the presheaf $U \rar
\widehat{\text{bar}^{\bullet}(\dif(\cale)(U))}$. The terms of this
complex are the same as in $\widehat{\text{hoch}}(\dif(\cale))$ but
the differential is the bar differential. We denote this by $\widehat{\text{bar}}(\dif(\cale))$. \\

We can obtain  the completed Tsygan's double complex of
$\dif(\cale)$. This is a double complex of sheaves of
$\compl$-vector spaces on $X$, whose non-positive even columns are $
\widehat{\text{hoch}}^{\bullet}(\dif(\cale))$ and whose negative odd
columns are $\widehat{\text{bar}}(\dif(\cale))$. The horizontal
differentials in Tsygan's double complex extend to this situation to
give us well defined horizontal differentials. We will denote the
total complex of this
double complex by $\widehat{\text{Cycl}}^{\bullet}(\dif(\cale))$. \\

Before we proceed, we note that as before, we have a (degree
preserving) natural map of complexes $\gamma:
\widehat{\text{Cycl}}^{\bullet}(\dif(\cale)) \rar
\widetilde{\text{Cycl}}^{\bullet}(\dif(\cale))$. This induces a map
$\gamma_*:
\hcoh{\bullet}{\widehat{\text{Cycl}}^{\bullet}(\dif(\cale))} \rar
\hcoh{\bullet}{\widetilde{\text{Cycl}}^{\bullet}(\dif(\cale))}$.
Note that the map $S$ extends to a map of complexes
$\widehat{\text{Cycl}}^{\bullet}(\dif(\cale)) \rar
\widehat{\text{Cycl}}^{\bullet +2}(\dif(\cale))$. Similarly, the map
${\calg I}$ extends to a map of complexes ${\calg I}:
\widehat{\text{hoch}}(\dif(\cale)) \rar
\widehat{\text{Cycl}}^{\bullet}(\dif(\cale))$
.Further, we have the following commutative diagrams. \\

$$\begin{CD}
 \widehat{\text{hoch}}(\dif(\cale)) @>{\calg I}>> \widehat{\text{Cycl}}^{\bullet }(\dif(\cale)) \\
 @V{\gamma}VV    @VV{\gamma}V \\
 \widetilde{\text{hoch}}(\dif(\cale)) @>{\calg I}>> \widetilde{\text{Cycl}}^{\bullet}(\dif(\cale))\\
 \end{CD} $$

 $$\begin{CD}
 \widehat{\text{Cycl}}^{\bullet}(\dif(\cale)) @>{S}>> \widehat{\text{Cycl}}^{\bullet + 2}(\dif(\cale)) \\
 @V{\gamma}VV    @VV{\gamma}V \\
 \widetilde{\text{Cycl}}^{\bullet}(\dif(\cale)) @>{S}>> \widetilde{\text{Cycl}}^{\bullet+2}(\dif(\cale))\\
 \end{CD} $$

We denote by $\ccccyc{-j}{(\dif(\cale))}$ the hypercohomology $\hcoh{-j}{\widehat{\text{Cycl}}^{\bullet}(\dif(\cale))}$. \\

Again we note that the complex
$\widetilde{\text{Cycl}}^{\bullet}(\dif(\cale))$ is a complex of
$\underline{\text{C}}^{\infty}$ modules. Thus, its hypercohomology
is precisely $\ccyc{\bullet}{(\dif(\cale))}$. It follows that we
have the following corollaries of Propositions 13 and 14
respectively.

\begin{cor} (Corollary to Proposition 13)
There exist traces $\hat{\text{tr}}_{2i}: \ccccyc{-2i}{(\dif(\cale))} \rar \compl$. \\
\end{cor}

\begin{proof} $\hat{\text{tr}_{2i}} = \tilde{\text{tr}_{2i}}
\circ \gamma_* $.
\end{proof}

\begin{cor} (Corollary to Proposition 14)
The following diagrams commute.\\

$$\begin{CD}
\choch{0}{(\dif(\cale))} @>{\calg I}>> \ccccyc{0}{(\dif(\cale))} \\
@VV{\hat{\text{tr}}}V @V{\hat{\text{tr}}_{0}}VV \\
\compl @> \id >> \compl \\
\end{CD} $$

$$\begin{CD}
\ccccyc{-2i}{(\dif(\cale))} @>{S}>> \ccccyc{-2i+2}{(\dif(\cale))} \\
@VV{\hat{\text{tr}}_{2i}}V @V{\hat{\text{tr}}_{2i-2}}VV \\
\compl @> \id >> \compl \\
\end{CD} $$
\end{cor}

\begin{proof} This is an immediate consequence of Proposition 14 and the commutative diagrams shown immediately before Corollary 7 . \\
\end{proof}

We have the following cyclic homology analog of Lemma 3. \\

\begin{lem} $\widehat{\text{Cycl}}^{\bullet}(\dif(\cale))$ is quasiisomorphic to $\underline{\compl}[2n] \oplus \underline{\compl}[2n+2] \oplus .... $
\end{lem}

\begin{proof} It suffices to verify this assertion locally at the level of pre sheaves. In other words, we need to look at the
completed Tsygan's double
complex for $\dif(\cale)(U)$ where $U$ is an open set so that
$\cale$ is trivial on $U$. Of course, all tensor products here are
completed tensor products. \\

The completed bar complex of $\dif(\cale)(U)$ is acyclic since the
homotopy of recollection 1, Section 2.1 between the identity
endomorphism and the $0$ endomorphism of the bar complex of
$\dif(\cale)(U)$ can be shown to "extend to" a homotopy between the
identity and $0$ endomorphisms of the completed bar complex of
$\dif(\cale)(U)$. The complex
$\widehat{\text{hoch}(\dif(\cale)(U))}$ is quasiisomorphic to
$\compl[2n]$ by Lemma 3. It follows that the spectral sequence
converging to the completed cyclic homology of $\dif(\cale)(U)$ that
arises out of the filtration of the completed
 Tsygan's double complex for $\dif(\cale)(U)$ by columns satisfies $$E_1^{-2i,-2n} \simeq \compl \text{ } \forall i \geq 0$$
 $$E_1^{p,q} = 0 \text{  otherwise. } $$ The desired lemma follows immediately from this.    \\

\end{proof}

\begin{cor} $\ccccyc{-2i}{(\dif(\cale))} = \hh{2n}{\compl} \oplus \hh{2n-2}{\compl} \oplus... \oplus
\hh{2n-2i}{\compl} $.
\end{cor}

Denote the isomorphism in the above corollary by $J_{2i}$.  Then, following
the proof of lemma 5, we obtain the following proposition. \\

\begin{prop}
The following diagrams commute. \\
$$\begin{CD}
\ccccyc{-2i}{(\dif(\cale))} @>>{J_{2i}}> \hh{2n}{\compl} \oplus
\hh{2n-2}{\compl} \oplus... \oplus \hh{2n-2i}{\compl} \\
 @VV{S}V             @VVV\\
 \ccccyc{-2i+2}{(\dif(\cale))} @>>{J_{2i-2}}>  \hh{2n}{\compl} \oplus \hh{2n-2}{\compl}
\oplus... \oplus \hh{2n-2i+2}{\compl} \\
\end{CD} $$

$$\begin{CD}
\choch{0}(\dif(\cale)) @>{\calg I}>> \ccccyc{0}(\dif(\cale))\\
@VV{\beta_{\cale}}V     @VV{J_0}V \\
 \hh{2n}{\compl}     @>id>> \hh{2n}{\compl} \\
 \end{CD} $$

The vertical arrow in the right column of the first diagram is the obvious projection. \\
\end{prop}

\subsection{Proof of theorem 3}

\begin{proof}(Proof of Theorem 3)

We apply the second diagram in Corollary 8 $i$ times to obtain the following diagram. \\

$$\begin{CD}
\ccccyc{-2i}{(\dif(\cale))} @> {S^i}>> \ccccyc{0}{(\dif(\cale))} \\
@VV{\hat{\text{tr}}_{2i}}V      @V{\hat{\text{tr}}_{0}}VV \\
\compl  @> \id >> \compl \\
\end{CD} $$

We apply the first diagram Proposition 15 $i$ times to obtain the following diagram. \\

$$\begin{CD}
\ccccyc{-2i}{(\dif(\cale))} @> {S^i}>> \ccccyc{0}{(\dif(\cale))} \\
@VV{J_{2i}}V                                                                          @V{J_{0}}VV \\
\hh{2n}{\compl} \oplus \hh{2n-2}{\compl} \oplus... \oplus \hh{2n-2i}{\compl} @>>> \hh{2n}{\compl} \\
\end{CD} $$

In the notation used to state Theorem 3, this tells us that
$I_{\cale, 2i,0} = I_{\cale ,0,0}$ and $I_{\cale, 2i, 2k} = 0$ for
$k > 0$. Lastly, the second diagram of Proposition 15 and the first
diagram of corollary 8 together imply that $\hat{\text{tr}} \circ
\beta_{\cale}^{-1} = \hat{\text{tr}_0} \circ J_0^{-1}$. This shows
that $I_{\cale,0,0} = I_{\cale}$, thereby completing the proof of
Theorem 3.\\

\end{proof}

Remark 1: Our notion of completed Hochschild and cyclic homologies
seems to differ in detail from that used in [FLS]. This forced us to
rework
some steps of [FLS], in particular, Proposition 6 and related matters, in our situation. \\

Remark 2: The key step in the proof of Theorem 1 in [FLS] consists
of showing on one hand that the linear functional defined in [FLS]
applied to the operator $\id$ in $\diff(\cale)$ is indeed its
supertrace i.e, the Euler characteristic of $\cale$, and showing
that the image of $\id$ in $\choch{0}{(\dif(\cale)}$ gives the class
${\text{ch}(\cale). \text{td}(T_X)}_{2n}$ ( this part is done by
citing [NT1] and [NT2]) after passing to $\hh{2n}{\compl}$ and then
applying the Hirzebruch R-R theorem.  In spite of the differences in
detail between our construction of completed Hochschild homology and
that of [FLS] it may
be checked that these two key steps go through, thus maintaining Theorem 1 in this situation. \\

\end{document}